\documentclass{hld2025} 
\def\R{{\mathbb{R}}}
\def\vx{{\mathbf{x}}}
\def\vy{{\mathbf{y}}}
\def\vz{{\mathbf{z}}}
\def\vA{{\mathbf{A}}}
\def\vB{{\mathbf{B}}}
\def\vC{{\mathbf{C}}}
\def\vT{{\mathbf{T}}}
\def\vX{{\mathbf{X}}}
\def\vY{{\mathbf{Y}}}
\def\mA{{\mathbf{A}}}
\def\mU{{\mathbf{U}}}
\def\mG{{\mathbf{G}}}
\def\mLambda{{\mathbf{\Lambda}}}
\def\R{{\mathbb{R}}}
\def\vomega{{\boldsymbol{\omega}}}
\def\dvom{{\dot{\boldsymbol{\omega}}}}
\usepackage{enumitem}
\usepackage{subcaption}

\title[Laplace Transform Viewpoint of Lookahead]{Understanding Lookahead Dynamics \\Through Laplace Transform}

\usepackage{xcolor}         
\definecolor{DarkGreen}{rgb}{0.1,0.5,0.1}
\definecolor{DarkRed}{rgb}{0.5,0.1,0.1}
\definecolor{DarkBlue}{rgb}{0.1,0.1,0.5}
\usepackage{hyperref}       
\hypersetup{
    unicode=false,          
    pdftoolbar=true,        
    pdfmenubar=true,        
    pdffitwindow=false,     
    pdfnewwindow=true,      
    colorlinks=true,       
    linkcolor=DarkBlue,          
    citecolor=DarkGreen,        
    filecolor=DarkRed,    
    urlcolor=DarkBlue,          
    pdftitle={},
    pdfauthor={},
}

\hldauthor{
\Name{Aniket Sanyal} 
\Email{aniket.sanyal@tum.de}\\
\addr {Technical University of Munich} 
\AND
\Name{Tatjana Chavdarova} \Email{tatjana.chavdarova@polimi.it}\\
\addr {Politecnico di Milano} }

\begin{document}

\maketitle

\begin{abstract}
We introduce a frequency-domain framework for convergence analysis of hyperparameters in game optimization, leveraging \emph{High-Resolution Differential Equations} (HRDEs) and Laplace transforms. Focusing on the Lookahead algorithm--characterized by gradient steps $k$ and averaging coefficient $\alpha$--we transform the discrete-time oscillatory dynamics of bilinear games into the frequency domain to derive precise convergence criteria. Our higher-precision $\mathcal{O}(\gamma^2)$-HRDE models yield tighter criteria, while our first-order $\mathcal{O}(\gamma)$-HRDE models offer practical guidance by prioritizing actionable hyperparameter tuning over complex closed-form solutions. Empirical validation in discrete-time settings demonstrates the effectiveness of our approach, which may further extend to locally linear operators, offering a scalable framework for selecting hyperparameters for learning in games.
\end{abstract}

\section{Introduction}

Saddle point optimization has become a major focus in machine learning due to its role in applications such as generative adversarial networks~\citep[GANs,][]{article}, adversarial training~\cite{goodfellow2015explainingharnessingadversarialexamples,madry2019deeplearningmodelsresistant}, and multi-agent reinforcement learning~\citep{shapley1953,littman1994,bertsekas2021rollout}. These are modeled as multiplayer games that aim to find a Nash equilibrium---a state where no player benefits from changing their strategy.

Traditional optimizers like Gradient Descent (GD) fail in these settings due to rotational dynamics inherent to games. Specialized algorithms~\citep[\textit{e.g.}, ][]{balduzzi2018mechanicsnplayerdifferentiablegames} may excel in purely rotational games but underperform in primarily potential games such as minimization. To bridge this gap, Lookahead was proposed as a wrapper that averages a base optimizer’s updates over $k$ steps, in a computationally efficient way; refer to \S~\ref{sec:prelim}. This enables divergent algorithms like GD to converge in saddle point problems.

Lookahead introduces two hyperparameters: the averaging frequency $k$ and the interpolation weight $\alpha\in(0,1)$. While most $(\alpha,k)$ combinations outperform the starting base optimizer, optimal joint tuning remains challenging. 
Existing methods optimize one parameter at a time, which works only for simple cases~\citep{ha2022lookahead}.
Understanding how to select these parameters jointly remains an open problem.

We analyze Lookahead’s dynamics relying on its \emph{High-Resolution Differential Equation}~\citep[HRDE,][]{chavdarova2023hrdes} representation, and using frequency domain tools.
Our approach is based on the following insight:
\textit{(i)} while the joint iterate cycles around the solution due to the rotational game dynamics,
\textit{(ii)} each player individually exhibits periodic behavior;
see Figure~\ref{fig:gd_ilustration} for intuition.
This paper studies:
\begin{center}
\emph{
Can estimating player-specific oscillation frequencies guide joint selection of Lookahead’s hyperparameters? }
\end{center}

To address this, we derive a general HRDE for Lookahead with arbitrary $\alpha$ and $k$, extending prior work limited to $k=2,3$~\citep{chavdarova2023hrdes}. Retaining up to linear terms in the step-size $\gamma$, $\mathcal{O}(\gamma)$-HRDE, we formulate a unified HRDE that captures the interplay between  $\alpha$ and $k$, and rotational dynamics.

We then analyze the general-$(\alpha,k)$  $\mathcal{O}(\gamma)$--HRDE using Laplace transforms, a tool that converts differential equations into algebraic forms in the frequency domain. This ``frequency dual'' representation simplifies stability analysis and reveals how oscillations (frequency/amplitude) depend on $(\alpha,k)$, guiding hyperparameter selection.

\paragraph{Contributions.} Our contributions can be summarized as follows.

\begin{itemize}[leftmargin=*,itemsep=0cm,parsep=0.05cm]
    \item We derive Lookahead's general-$(\alpha,k)$  $\mathcal{O}(\gamma)$--HRDE. 

    Using the Laplace transform, we obtain the frequency-domain representation of the $\mathcal{O}(\gamma)$--HRDE and discuss its interpretation.

    \item For two-player bilinear games---see \eqref{eq:bilinear_game}---we derive a closed-form trajectory of Lookahead’s dynamics over time. By analyzing the poles of its frequency-domain representation, we establish parameter constraints (on $\alpha,k$, and the step size) for convergence and empirically validate these in discrete-time implementations on high rotation games; see \eqref{eq:quadratic_game}. However, in games with weaker rotations (due to stronger potential terms, e.g., quadratic components), convergence may occur even when our initial conditions are violated.

    \item To understand this disparity, we show that refining the analysis using an even higher-resolution HRDEs (deriving an $O(\gamma^2)$-HRDE) yields ``tighter'' convergence criteria for primarily potential games. In particular, we derive the corresponding convergence condition and empirically validate it relative to the $O(\gamma)$-HRDE one. 

\end{itemize}

The proposed approach is general and can be applied and extended to other optimization algorithms that can be represented using a (high-resolution) ordinary differential equation.

\begin{center}
\noindent
\begin{minipage}[t]{0.59\linewidth}
\label{fig:gd_ilustration}
    \includegraphics[width=\linewidth,trim= 5em .5em .3em .1em,clip]{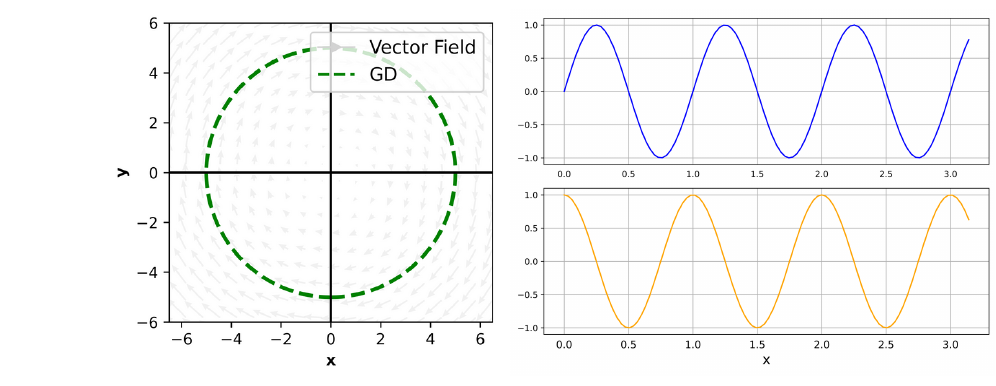}
\end{minipage}
\hfill
\begin{minipage}[t]{0.38\linewidth}
\vspace{10pt}
    \textbf{Figure 1:} \textbf{GD on game~\eqref{eq:bilinear_game} (infinitesimal $\gamma$).}
    In the joint parameter space, GD oscillates circularly around the solution.
    In their own parameter space, each player oscillates harmonically with a right-angled phase difference between them (right).
\end{minipage}
\end{center}

\paragraph{Related works.}
 \citet{Anagnostides2021FrequencyDomainRO} introduced a frequency-domain framework to analyze the convergence of algorithms that incorporate historical information, such as past gradients. They employed bilateral $Z$-transform to study the convergence of algorithms such as the \emph{Optimistic Gradient Descent}~\citep{popov1980}. It simplifies convergence analysis and establishes linear convergence under the co-coercive and monotone operators. 
 Similar to their approach, which ensures convergence by restricting the poles to lie within the unit circle in the Z-plane, our work imposes conditions on the poles in the complex frequency domain instead.
 
\section{Prelimnaries}\label{sec:prelim}
Small letters denote variables in $\mathbb{R}$.
Small and capital bold letters denote vectors and matrices, respectively.
$*$ denotes the convolution operator defined in Appendix~\ref{app:background}. 

\paragraph{Problem.}
Our study includes the following problem:
\begin{equation}
    \tag{QD}
    \label{eq:quadratic_game}
    \min_{\substack{\vx}} \max_{\substack{\vy}}  \quad 
    \langle \vx, \vA\vy \rangle  + \frac{1}{2}\langle \vx, \vC\vx \rangle - \frac{1}{2} \langle \vy, \vB\vy \rangle 
\end{equation}
where $\vx, \vy \in \mathbb{R}^{\frac{d}{2}}$ and $\vA, \vB, \vC $ are $\frac{d}{2} \times \frac{d}{2}$ positive semi-definite matrices, and full rank unless otherwise specified.

We will also consider a special case of \eqref{eq:quadratic_game} as follows:
\begin{equation}
    \tag{BG}
    \min_{\substack{ \vx \in \mathbb{R}^{\frac{d}{2}}}} \quad \max_{\substack{\vy \in \mathbb{R}^{\frac{d}{2}}}} \quad \vx^{\intercal} \mA \vy \,.
\label{eq:bilinear_game}
\end{equation}

\eqref{eq:quadratic_game} and \eqref{eq:bilinear_game} are an instance of Variational Inequality~\citep[VI,][]{stampacchia_formes_1964,facchineiFiniteDimensionalVariationalInequalities2003} problems defined by a tuple of the operator $F\colon \mathbb{R}^n \to \mathbb{R}^n$ and its domain---herein $\R^n$---as follows:
\begin{equation} \label{eq:vi} \tag{VI}
  \emph{find } z^\star \quad\text{s.t.}\quad \langle \vz-\vz^\star, F(\vz^\star)\rangle \geq 0, \quad \forall \vz \in \mathbb{R}^n \,.
\end{equation}

For \ref{eq:bilinear_game}, we can compactly write the operator as $F \equiv 
(\nabla_x (\vx^{\intercal} A \vy), \ \nabla_y (-\vx ^{\intercal} A\vy) )^\intercal = (A\vy, -A\vx)^\intercal \,.$

\paragraph{Optimization methods.}
\emph{Gradient descent} with step size $\gamma \in [0,1]$  for VIs is then:
\begin{equation} \tag{GD}\label{eq:gd}
    \vz_{n+1} = \vz_n - \gamma F(\vz_n) \,. 
\end{equation}
The \emph{lookahead} algorithm~\citep{zhang2019lookaheadoptimizerksteps,chavdarova2021lamm} computes its update as a point on the line between the current iterate ($\vz_n$) and the iterate obtained after taking $k \geq 1$ steps with the base optimiser ($\tilde \vz_{t+k}$): 
\begin{equation}
\tag{LA}
\vz_{n+1} \leftarrow \vz_n + \alpha (\tilde  \vz_{n+k} -  \vz_n), \quad \alpha \in [0,1]  
\,. \label{eq:lookahead}
\end{equation}
\noindent
In addition to the step size of the base optimizer (herein \ref{eq:gd}), Lookahead has two hyperparameters: 
\begin{itemize}
\item $k$---the number of steps of the base optimiser used to obtain the prediction $\smash{\tilde x_{t+k}}$; and
\item $\alpha$---averaging parameter that controls the step toward the predicted iterate  $\smash{\tilde x}$. The larger the closest, and $\alpha=1$~ \eqref{eq:lookahead} is equivalent to running the base algorithm only (has no impact). 
\end{itemize}
We provide the general pseudocode for Lookahead in Appendix \ref{app:pseudocode}.
\paragraph{Linear system stability.}
A linear dynamical system  $(\dot{\vz}, \dvom)^\intercal = \mathbf{C} (\vz, \dvom)$ 
is considered stable if and only if all eigenvalues of \( C \) have negative real parts, i.e., \( \Re(\lambda_i) < 0 \) for every \( \lambda_i \in \text{Sp}(\mathbf{C}) \). Additionally, the \emph{Routh–Hurwitz} stability criterion provides a necessary and sufficient condition for determining the stability of such a system. Rather than computing the eigenvalues explicitly, this criterion evaluates stability using the coefficients of the characteristic polynomial by constructing the Hurwitz array or its generalized form~\citep{xie1985} when dealing with complex coefficients.

\paragraph{Laplace transform.}
Converts a function that depends on a real variable (the time $t$) to a function with a complex variable---the complex-valued frequency $s$. 
The Laplace transform of $\vx(t)$ is: 
\begin{equation}
\label{eq:Laplace}\tag{LT}
\mathcal{L}\{\vx(t)\} = \mathbf{X}(s) = \int_{0}^{\infty} \vx(t) \ e^{-st} \, dt \,, 
\end{equation}  
where $s = \sigma + j \omega$ is the complex frequency. 
$\sigma$ is the real exponential decay factor, while $\omega$ represents the frequency. 

Unlike the Fourier transform, the Laplace transform includes $\sigma$, enabling the analysis of exponentially growing or unbounded systems.
By converting differential equations to algebraic forms, it simplifies the solution derivation. These solutions are inverted back to the time domain using the \emph{inverse Laplace transform} (see Appendix~\ref{app:background}).

\section{General-$(\alpha,k)$  $\mathcal{O}(\gamma)$--HRDE of Lookahead \& Convergence on \ref{eq:bilinear_game}
}\label{sec:la_general_hrde}

This section presents the Lookahead HRDE with $(\alpha, k)$ kept as parameters. We defer its derivation to Appendix \ref{app:la_hrde_derivation}. 
The $\mathcal{O}(\gamma)$-HRDE for~\eqref{eq:lookahead} with \eqref{eq:gd} as a base optimizer is:
\begin{equation}
\tag{LA-$\gamma$-HRDE}
    \ddot{\vz}(t)     = - \frac{2}{\gamma}\cdot \dot{\vz}(t) - \frac{2k\alpha}{\gamma} \cdot F(\vz(t)) + \big( \sum_{i=1}^{k-1} i \big) \cdot 2 \alpha \cdot \mathcal{J}(\vz(t)) \cdot F(\vz(t)) \,,
\label{eq:la-hrde} 
\end{equation}
where $\vz \in \mathbb{R}^n$ is the vector of all players, 
$k$ is the number of \eqref{eq:gd} steps  before the $\alpha$-averaging, $\gamma$ is the \eqref{eq:gd} step size, $F(\cdot)$ is the operator and $\mathcal{J}(\cdot)$ is the Jacobian, that is, the derivative of $F$.

 \begin{theorem}{(Convergence on bilinear games)}
\label{theo:convergence}
    For any \(\gamma\), the real part of the eigenvalues of $\mathbf{C}_{LAk-GDA}$ is always negative, i.e., \( \Re(\lambda_i) < 0, \forall \lambda_i \in \text{Sp}(\mathbf{C}_{LAk-GDA}) \).  
Thus, the \( O(\gamma) \)-HRDE of the LA k-GDA method (\ref{eq:la-hrde}) converges on the BG problem for any step size.
\end{theorem}

\noindent
The proof is provided in Appendix \ref{sec:theorem1-proof}. While Theorem \ref{theo:convergence} guarantees convergence, it offers no guidance for selecting hyperparameters or their relationship to the step size. We address this gap below.

\section{Lookahead for \ref{eq:bilinear_game}: Solution Trajectory \& Convergence}

We compute \eqref{eq:la-hrde} frequency dual via Laplace transform \eqref{eq:Laplace}, under the bilinear model \eqref{eq:bilinear_game}, then derive the solution trajectory with its inverse Laplace transform.

\begin{theorem}[Solution trajectories]
\label{theo:2}
Consider the~\ref{eq:bilinear_game} (with matrix $\mA$). Let $U$ be the orthogonal matrix from eigen decomposition of $\mA$ i.e., $\mA = \mU \mLambda \mU^\intercal$, where $\mLambda = diag(\lambda_i)$.
    The trajectory of the individual players $\vx(t)$ and $\vy(t)$ of the~\eqref{eq:la-hrde} continuous time dynamic with parameters $(k, \alpha)$ and Gradient Descent with step size $\gamma$ as the base optimizer, are as follows:
    \begin{equation}\tag{$x$-Sol}
    \vx(t) = -\frac{2k\alpha}{\gamma}  (\mG \ast \vy)(t) + \mU \, \mathrm{diag}\left( 
  e^{-\frac{t}{\gamma}} \left[ 
  \cosh(\omega_i t) \, x_i(0) + \frac{1}{\omega_i} \sinh(\omega_i t) \big( \dot{x}_i(0) + \frac{1}{\gamma} x_i(0)
  \big)
  \right]\right) \mU^\intercal \,,
    \label{eq:trajectory_1}
    \end{equation}
    \begin{equation}\tag{$y$-Sol}
       \vy(t) =  \frac{2k\alpha}{\gamma} (\mG * \vx)(t) + \mU \, \mathrm{diag}\left(
  e^{-\frac{t}{\gamma}} \left[ 
  \cosh(\omega_i t) \, y_i(0) + \frac{1}{\omega_i} \sinh(\omega_i t) \big( \dot{y}_i(0) + \frac{1}{\gamma} y_i(0)
  \big)
  \right]\right) \mU^\intercal \,,
    \end{equation}
where $*$ denotes the convolution operator \eqref{convolution},  $x_i$ is the $i$th component of $\vx$,
$\omega_i = \sqrt{\frac{1}{\gamma^2} - \alpha k(k-1) \lambda_i}$, $\mG(t) = \mU \, \mathrm{diag}\left( \frac{1}{\omega_i} e^{-\frac{t}{\gamma}} \sinh(\omega_i t) \right) \mU^\intercal$,  $(\vx(0), \vy(0))$  are the initial points and $\dot{\vx(0)}, \dot{\vy(0)}$ are the initial momentum, respectively.
\end{theorem}

\noindent\textbf{Interpretation.} 
Exponential terms ensure convergence by reducing oscillation amplitude over time; periodic terms capture oscillatory behavior around the equilibrium--a hallmark of minimax optimization. These dynamics are mutually coupled: $ \vx(t)$ depends on $\vy(t)$ and vice versa, with oscillation frequency governed by $\omega_0\,.$

For \eqref{eq:bilinear_game}, we derive a decoupled trajectory in Appendix \ref{app:la_bg_convergence} where $ \vx(t)$ becomes independent of $\vy(t)$, reducing the problem to pure minimization over $\vx$. 

\subsection{GD \& Lookahead Convergence Analysis} \label{sec:convergence}

\paragraph{Gradient Descent~\eqref{eq:gd} for \ref{eq:bilinear_game}.}

Through Laplace transform of \eqref{eq:gda-hrde}, we move to the frequency domain and enforce pole constraints to study convergence; see Appendix \ref{app:gda_convergence_analysis}. In line with existing results~\citep{chavdarova2023hrdes}, our analysis shows divergence of \eqref{eq:gda-hrde} for all $\gamma > 0$.

\paragraph{Lookahead~\eqref{eq:lookahead} for \ref{eq:bilinear_game}.}
We transform \ref{eq:la-hrde} in the complex frequency domain using \eqref{eq:Laplace}. This gives separate equations for the players, which are then substituted (see Appendix \ref{app:la_bg_convergence}): 
\begin{equation}
    \label{eq:freq_hrde} \tag{$X$--TF}
    \mathbf{X}(s) = \dfrac{-2k\gamma\alpha \vy(0) (s+ \frac{2}{\gamma})\vA -2k\alpha\gamma\dot{\vy}(0)\vA}{\gamma^2 (s^2 +\frac{2s}{\gamma}+\alpha k(k-1) \vA^2)^2 + 4k^2\alpha^2 } + \dfrac{(\dot{\vx}(0) + (s+\frac{2}{\gamma})\vx(0)) \gamma^2(s^2 +\frac{2s}{\gamma}+\alpha k(k-1))}{\gamma^2 (s^2 +\frac{2s}{\gamma}+\alpha k(k-1)\vA^2)^2 + 4k^2\alpha^2} \,, 
\end{equation}
for the $\vx$-player; refer to \ref{app:la_bg_convergence} for the $\vy$-player expression.

\noindent\emph{Interpretation.}
Notably, for \ref{eq:bilinear_game}, substituting $\mathbf{Y}(s)$ in the transfer function $\mathbf{X}(s)$, eliminates interdependence on $\mathbf{Y}(s)$. $\mathbf{X}(s)$ and $\mathbf{Y}(s)$ depend only on the opponent's initial values. Thus, a time-dependent trajectory function $\vx(t)$ exists---through the inverse Laplace transform---that does \emph{not} depend on the trajectory of the opponent $\vy(t)$. This observation extends to affine VIs; see  Appendix~\ref{app:la_bg_convergence}

\noindent

Importantly, \eqref{eq:Laplace} provides convergence conditions without requiring the initial values or solution trajectory of the corresponding HRDEs.
Convergence is determined by analyzing the poles of $\mathbf{X}(s)$. If the dominant pole (largest real part) lies in the negative half-plane, the trajectory $\vx(t)$ converges to $0$ as $t \to \infty$. For Lookahead, convergence for \ref{eq:bilinear_game} requires:

\begin{equation}\label{eq:cond_hamiltonian}\tag{BG-Cond}
    \alpha < \frac{k - 1}{k} \,. 
\end{equation}

\paragraph{Lookahead~\eqref{eq:lookahead} for \ref{eq:quadratic_game}.}
Using the above approach, in Appendix \ref{app:la_bg_convergence} we derive the LA convergence condition for \ref{eq:quadratic_game} (refer to \eqref{eq:convergence-inequality-full}). Furthermore, it can be seen therein that the latter reduces to the former~\ref{eq:cond_hamiltonian} as a special case. 

\setcounter{figure}{1} 
\begin{center}
\begin{minipage}[t]{0.65\linewidth}
    \includegraphics[width=\linewidth]{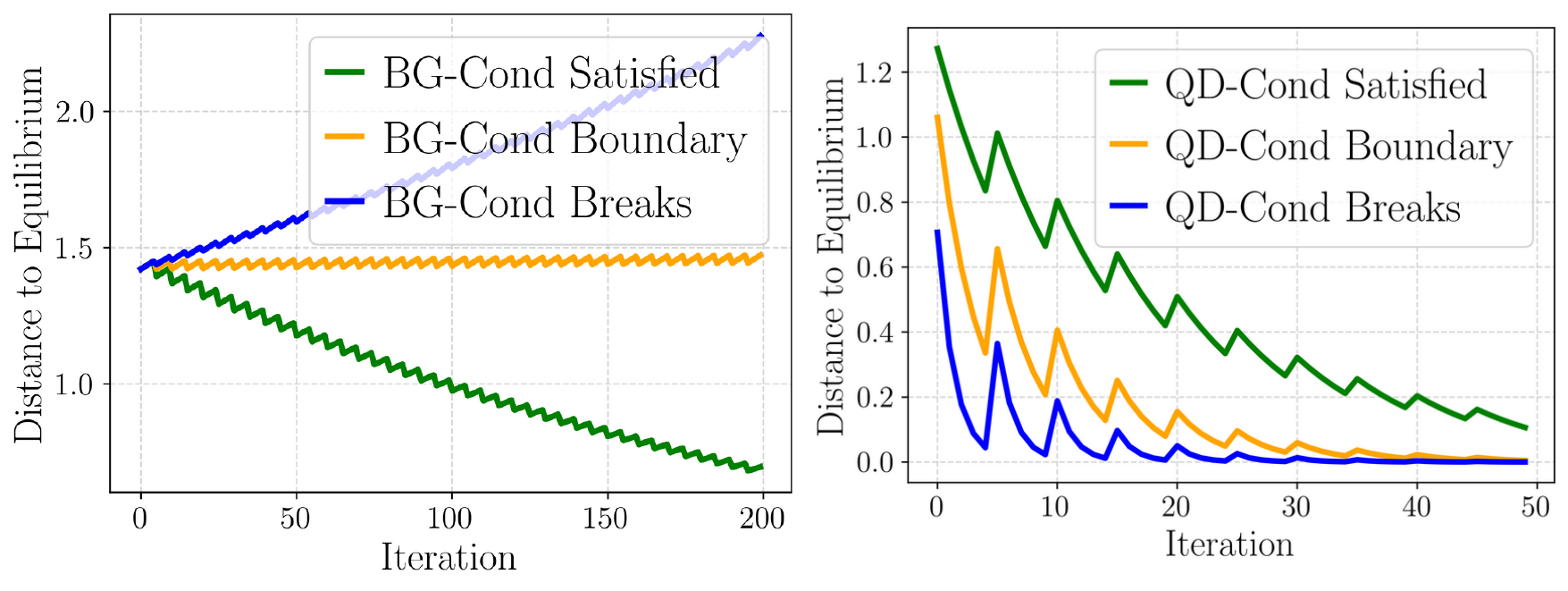}
\end{minipage}
\hfill
\begin{minipage}[t]{0.33\linewidth}
     \refstepcounter{figure}
     \vspace{-100pt}
    \textbf{Figure 2:} Verification of the obtained convergence conditions on: \textbf{left}--\ref{eq:convergence-condition} yields a tight bound for \ref{eq:bilinear_game}; \textbf{right}--\ref{eq:convergence-inequality-full} is not a tight bound for \ref{eq:quadratic_game}.
    See \S~\ref{sec:convergence} and Appendix \ref{app:experiments-1}.
    \label{fig:convergence-divergence}
\end{minipage}
\end{center}
 
\noindent
\textbf{Experiments summary--Figure~\ref{fig:convergence-divergence}.}
Transforming to the frequency domain yields a tighter convergence bound for \eqref{eq:bilinear_game}; violating it causes divergence. However, in highly potential games (as $\vB$ and $\vC$ grow less sparse), convergence may persist despite violations. Our analysis shows that deriving a higher-order HRDE condition also predicts divergence for such games. Yet, since divergence is more critical, even a weaker (but more tractable) condition retains practical utility.

\section{\(O(\gamma^2)\)--HRDE for Lookahead \& Convergence}
\label{sec:ogamma-2}
In Appendix \ref{app: hrde-derive-second} we derive the general $(\alpha,k)$  $\mathcal{O}(\gamma^2)$--HRDE, which relative to the \ref{eq:la-hrde} has the following additional terms on the right hand side: $(\sum_{i=0}^{k-1} i^2) \gamma F(\vz(t))^\intercal H(\vz(t))F(\vz(t)) - 2(k-2)^2\gamma J(\vz(t))^2 F(\vz(t))\,,$ where $H$ is the derivative of the Jacobian $J$; full expression in \eqref{eq:LAk-HRDE-second}.

\noindent\textbf{Convergence \& convergence condition of \ref{eq:LAk-HRDE-second} for a special case of \eqref{eq:quadratic_game}.}
In Appendix ~\S \ref{sec:Convergence_Analysis_ogamma2} we prove convergence of \ref{eq:LAk-HRDE-second} for the quadratic game defined in Figure~\ref{fig:3}.
Applying the same pole analysis and Routh-Hurwitz criterion, we derive condition \ref{eq:convergence-o-gamma2}. Due to the complexity and lack of intuitive insight from the analytic form, we evaluate it numerically. Setting $\gamma^2 \rightarrow 0$ recovers \eqref{eq:convergence-inequality-full} with $\vB = \vC = (1-\beta) I$ and $\vA = \beta I$.

\noindent
\textbf{Experiments summary--Figure~\ref{fig:3}.}
We find the \(O(\gamma^2)\)-HRDE convergence condition more accurate for predominantly potential games, enabling optimal tuning of \(\gamma\) with fixed \(k, \alpha\). 
See Appendix \ref{app:experiments-2}.

\begin{center}
\begin{minipage}[t]{0.48\linewidth}
    \includegraphics[width=\linewidth,trim= 0em 0em 0em .5em,clip]{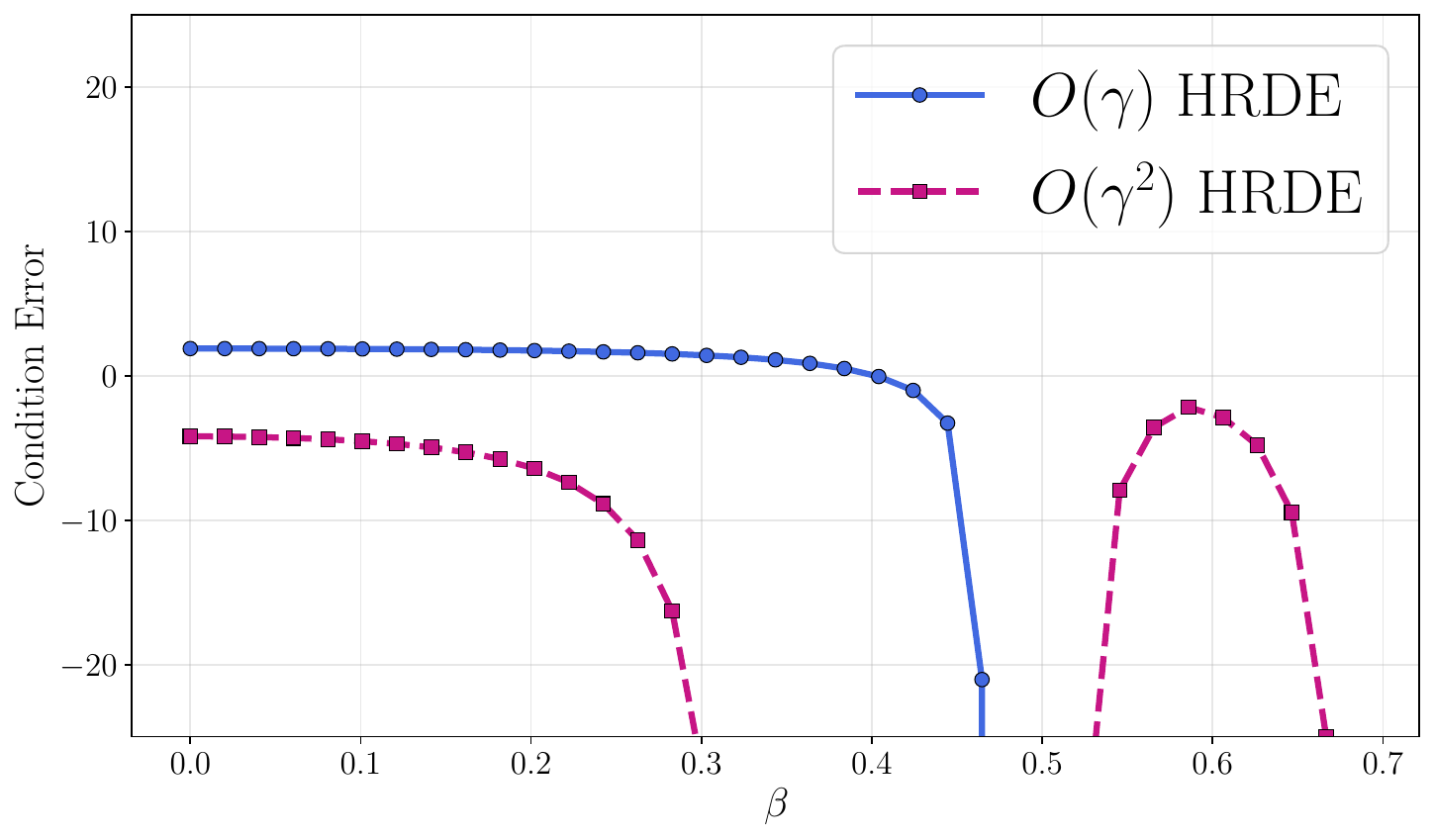}
\end{minipage}
\hspace{.5em}
\begin{minipage}[t]{0.48\linewidth}
\vspace{-122pt}
    \textbf{Figure 3:} \textit{Condition error for $\beta \in [0,1]$ of:}
    \begin{align*}
        \min_{\substack{x} \in \mathbb{R}} \max_{\substack{y} \in \mathbb{R}} \ (1-\beta) x^2 + \beta x y - (1-\beta)y^2 \,.   
    \end{align*}
    Negative error is preferred (Appendix \ref{app:experiments-2}); a positive error indicates subpar performance under convergence conditions. Hyperparameters tuned using our $O(\gamma^2)$ HRDE condition achieve negative error across the game spectrum.
    \refstepcounter{figure}
    \label{fig:3}
\end{minipage}
\end{center}

\section{Discussion}
Standard hyperparameter tuning in the joint parameter space struggles in game optimization; we address this with a frequency-domain framework using higher-precision ODEs to derive convergence conditions. By leveraging local linearity in complex dynamics, our approach enables practical Lookahead hyperparameter selection. Crucially, even simplified models (e.g., $\mathcal{O}(\gamma)$-HRDE) retain value: they prioritize preventing divergence--a critical practical concern--over analytically exhaustive precision, offering actionable guidance for stable training.

\section*{Acknowledgments}
TC was supported by the FAIR (Future Artificial Intelligence Research) project, funded by the NextGenerationEU program within the PNRR-PE-AI scheme (M4C2, Investment 1.3, Line on Artificial Intelligence).

\bibliography{main}
\clearpage
\appendix
\section{Additional Preliminaries}
In this section, we present the pseudocode for the two-player Lookahead algorithm introduced in~\S\ref{sec:prelim}, detailing its update steps and interaction between players. This serves as a concrete reference for the theoretical framework discussed in the main text. Additionally, we include supplementary background material—such as a review of the convolution operator—to provide the necessary analytical tools and intuition that support the derivations and interpretations in the main paper.
\subsection{Pseudocode for General Lookeahead}
\label{app:pseudocode}
\begin{algorithm2e}[H]
\caption{General Lookahead--Minmax Pseudocode}
\label{alg:lookahead_minimax}
\DontPrintSemicolon
\KwData{Stopping time $T$, learning rates $\eta_\theta, \eta_\phi$, initial weights $\theta_0, \phi_0$, lookahead $k$ and $\alpha$, base optimizer updates $\delta_{(\cdot)}^{\phi}, \delta_{(\cdot)}^{\theta}$, real data distribution $p_d$, noise--data distribution $p_z$}
\KwResult{Final weights $\theta_T$, $\phi_T$}
\BlankLine
$\tilde \theta_{0,0} \gets \theta_0$\;
$\tilde \phi_{0,0} \gets \phi_0$\;
\For{$t \gets 0$ \KwTo $T{-}1$}{
    \For{$i \gets 0$ \KwTo $k{-}1$}{
        Sample $\mathbf{x} \sim p_d$, $\mathbf{z} \sim p_z$\;
        $\tilde\phi_{t,i+1} \gets \tilde \phi_{t,i} - \eta_\phi \delta_{t,i}^{\phi} (\tilde \theta_{t,i}, \tilde\phi_{t,i}, \mathbf{x}, \mathbf{z})$\;
        Sample $\mathbf{z} \sim p_z$\;
        $\tilde \theta_{t,i+1} \gets \tilde \theta_{t,i} - \eta_\theta \delta_{t,i}^{\theta} (\tilde\theta_{t,i}, \tilde\phi_{t,i}, \mathbf{z})$\;
    }
    $\phi_{t+1} \gets \phi_{t} + \alpha_{\phi} (\tilde\phi_{t,k} - \phi_t)$\; \label{alg:backtrack_phi}
    $\theta_{t+1} \gets \theta_{t} + \alpha_{\theta} (\tilde\theta_{t,k} - \theta_t)$\; \label{alg:backtrack_theta}
    $\tilde\theta_{t+1,0} \gets \theta_{t+1}$\;
    $\tilde\phi_{t+1,0} \gets \phi_{t+1}$\;
}
\Return{$\theta_T$, $\phi_T$}\;
\end{algorithm2e}

\subsection{Extended Background}\label{app:background}

\paragraph{Convolution operator.}
The convolution operator is an important operator in control theory, for systems where the output at a given instance depends on the cummulative effect of past inputs.
Let \( f(t) \) and \( g(t) \) be two functions defined for \( t \geq 0 \), that is, \( f(t)\colon [0, \infty) \to \mathbb{R} \) and \( g(t)\colon [0, \infty) \to \mathbb{R} \). The convolution of \( f(t) \) and \( g(t) \) is given by the integral:
\begin{equation}
\label{convolution} \tag{CONV}
    (f * g)(t) = \int_{0}^{t} f(\tau) g(t - \tau) \, d\tau, \quad t \geq 0 \,.
\end{equation}

\paragraph{Laplace transform.}
For a function $\vx(t)$, the Laplace transform is defined as 
\begin{equation}
\tag{LT}
\mathcal{L}\{\vx(t)\} = \mathbf{X}(s) = \int_{0}^{\infty} \vx(t) e^{-st} \, dt
\end{equation}  
where $s$ represents the frequency and is a complex number and $t$ is the time.

The inverse Laplace transform is used to convert a function \( X(s) \) back into its original time-domain function \( x(t) \). It is defined by the Bromwich integral:

\begin{equation}
\label{inverse-laplace}\tag{iLT}
    \mathcal{L}^{-1}\{\mathbf{X}(s)\} = \vx(t) = \frac{1}{2\pi i} \int_{c-i\infty}^{c+i\infty} \mathbf{X}(s) e^{st} \, ds
\end{equation}

 The Laplace transform of the convolution of \( f(t) \) and \( g(t) \) is the product of their individual Laplace transforms. If \( \mathcal{L}\{f(t)\} = F(s) \) and \( \mathcal{L}\{g(t)\} = G(s) \), then:

\[
\mathcal{L}\{(f * g)(t)\} = F(s) G(s)
\]

This relation simplifies the convolution operation in the time domain to multiplication in the frequency domain.
The inverse Laplace transform of a product in the frequency domain yields the convolution in the time domain. That is,

\[
\mathcal{L}^{-1}\{F(s) G(s)\} = (f * g)(t) = \int_{0}^{t} f(\tau) g(t - \tau) \, d\tau
\]

Thus, the product of two functions in the frequency domain corresponds to the convolution of their inverse Laplace transforms in the time domain.

\section{\texorpdfstring{General $\mathcal{O}(\gamma)$--HRDE for Lookahead}{General O(gamma)--HRDE for Lookahead}}

\label{app:la_hrde_derivation}
In this section, we derive the $O(\gamma)$ HRDE for Lookahead as presented in ~\S \ref{sec:la_general_hrde}. This differential equation is a continuous time representation of the discrete Lookahead dynamics (\ref{eq:lookahead}). By transitioning to the continuous domain, we can leverage the well-established theory of differential equations to analyze the behavior of the algorithm and gain deeper insights into its underlying structure.

We follow the procedure of~\citet{chavdarova2023hrdes}, where we first re-write the methods in the following general form:
\begin{equation} \tag{GF} \label{eq:general_form}
  \frac{\vz_{n+1}-\vz_n}{\gamma} \!=\! \mathcal{U}(\vz_{n+k}, \dots, \vz_0) \,.
\end{equation}
We introduce the Ansatz $\vz_n \approx \vz(n \cdot \delta)$,  for some smooth
curve $\vz(t)$ defined for $t \ge 0$. A Taylor expansion gives:
\begin{align*}
    \vz_{n+1} &\!\approx\! \vz((n+1) \delta ) = \vz(n\delta) + \dot{\vz} (n\delta)\delta + \frac{1}{2} \ddot{\vz}(n \delta) \delta^2 + \dots \,.
\end{align*}
Thus,  when deriving $\mathcal{O}(\gamma)$--HRDEs, for the nominator of the left-hand side of \ref{eq:general_form} for all the methods considered in this work, we have:
\begin{equation} \tag{LHS-GF}\label{eq:time-taylor} 
    \vz_{n+1} - \vz_{n} \!\approx\!  \dot{\vz} (n\delta)\delta + \frac{1}{2} \ddot{\vz} (n \delta) \delta^2 + \mathcal{O}(\delta^3) \,.
\end{equation}

\paragraph{HRDEs for LA-GD}

The predicted iterates for $k=1, \dots, 4$, given \ref{eq:gd} as a base optimizer are as follows:
\begin{align*}
    \tilde\vz_{n+1} &= \vz_n - \gamma F(\vz_n)
    \tag{$\tilde\vz_{n+1}^{\text{GDA}}$}\label{eq:predicted_gda_1}\\
    \tilde\vz_{n+2} &= \tilde\vz_{n+1} - \gamma F(\tilde\vz_{n+1}) = \vz_n - \gamma F(\vz_n) - \gamma F(\vz_n - \gamma F (\vz_n) )
    \tag{$\tilde\vz_{n+2}^{\text{GDA}}$}\label{eq:predicted_gda_2}\\
    \tilde\vz_{n+3} &= \tilde\vz_{n+2} - \gamma F(\tilde\vz_{n+2}) \notag\\
    & = \vz_n - \gamma F(\vz_n) - \gamma F(\vz_n - \gamma F (\vz_n) ) 
    - \gamma F\big(
    \vz_n - \gamma F(\vz_n) - \gamma F(\vz_n - \gamma F (\vz_n) ) 
    \big)
    \tag{$\tilde\vz_{n+3}^{\text{GDA}}$}\label{eq:predicted_gda_3}\\
    \tilde\vz_{n+4} &= \tilde \vz_{n+3} - \gamma F(\tilde \vz_{n+3}) \notag\\
    &= \underbrace{
    \vz_n - \gamma F(\vz_n) - \gamma F(\vz_n - \gamma F (\vz_n) ) 
    - \gamma F\big(
    \vz_n - \gamma F(\vz_n) - \gamma F(\vz_n - \gamma F (\vz_n) ) 
    \big)
    }_{\tilde \vz_{n+3}} \notag\\
    & - \gamma F \Big(
    \underbrace{
    \vz_n - \gamma F(\vz_n) - \gamma F(\vz_n - \gamma F (\vz_n) ) 
    - \gamma F\big(
    \vz_n - \gamma F(\vz_n) - \gamma F(\vz_n - \gamma F (\vz_n) ) 
    \big)
    }_{\tilde \vz_{n+3}}
    \Big)
    \tag{$\tilde\vz_{n+4}^{\text{GDA}}$}\label{eq:predicted_gda_4}
\end{align*}

\paragraph{~\ref{eq:lookahead}2-GDA.}
The iterates of~\ref{eq:lookahead}2-GDA are obtained as follows:
\begin{align*}
    \vz_{n+1} = \vz_n + \alpha (\tilde\vz_{n+2} - \tilde\vz_n) = \vz_n + \alpha \Big(
    -\gamma F(\vz_n) - \gamma F\big(
    \vz_n - \gamma F(\vz_n)
    \big)
    \Big)\,.
\end{align*}
Using~\eqref{eq:time-taylor}, we get (where $\delta$ and $\gamma$ are the step sizes in time and parameter space, resp.):
\begin{align*}
    \frac{\dot \vz (n\delta) + \frac{1}{2} \delta^2 \ddot{\vz} (n) +\mathcal{O} (\delta^3) }{\gamma}
    = \alpha \Big( - 2F(\vz (n\delta)) + \gamma J(\vz(n\delta))F(\vz(n\delta))
    \Big) \,.
\end{align*}
Setting $\delta \!=\! \gamma$ and keeping the $\mathcal{O}(\gamma)$ terms yields:
\begin{align*}
\dot\vz(t) + \frac{\gamma}{2}\ddot\vz(t) = - 2\alpha F(\vz(t)) + \alpha \gamma J (\vz(t)) F(\vz(t)) \,,
\end{align*}

Writing it in phase-space representation we have:
\begin{equation}\tag{LA2-GDA-HRDE}\label{la2-gda-hrde}
\begin{aligned}
    \dot\vz(t) &= \vomega (t)\\
    \dot\vomega(t) &= -\frac{2}{\gamma} \vomega(t) - \frac{4\alpha}{\gamma} F(\vz(t)) + 2\alpha J(\vz(t)) F(\vz(t))
\end{aligned}
\end{equation}

\paragraph{~\ref{eq:lookahead}3-GDA.}
For~\ref{eq:lookahead}3-GDA using \eqref{eq:predicted_gda_3} we have:
\begin{equation} \label{eq:la3_gda_interm1}
\begin{aligned}
    \vz_{n+1} &= \vz_n + \alpha(\tilde \vz_{n+3} - \vz_n) \\
    &= \vz_n + \alpha\gamma\big[
    - F(\vz_n) - F\big(\vz_n - \gamma F(\vz_n)\big)
    - \underbrace{F\big( \vz_n - \gamma F(\vz_n) - \gamma F(\vz_n - \gamma F(\vz_n)) \big)}_{(\star)}
    \big]
\end{aligned}
\end{equation}
Similarly, by doing TE in coordinate space for ($\star$),  we get:
\begin{align*}
    F\big( \vz_n - \gamma F(\vz_n) - \gamma F(\vz_n - \gamma F(\vz_n)) \big) 
    &= F(\vz_n) - \gamma J (\vz_n)F(\vz_n) - \gamma J(\vz_n)F\big(\vz_n - \gamma F(\vz_n) \big) + \mathcal{O}(\gamma ^2) \\
    &= F(\vz_n) - 2\gamma J (\vz_n)F(\vz_n)  + \mathcal{O}(\gamma ^2) \,,
\end{align*}
wherein the second row we do additional TE of the last term in the preceding row.

Thus using \ref{eq:time-taylor} as well as replacing the above in~\eqref{eq:la3_gda_interm1} we have: 
\begin{align*}
\frac{\dot{\vz}(n\delta)\delta \!+\! \frac{1}{2} \ddot{\vz}(n\delta)\delta^2  + \mathcal{O}(\delta^3)}{\gamma} 
&=
\alpha \Big\{
- 3 F(\vz(n\delta)) + 3\gamma J (\vz(n\delta)) F(\vz(n\delta)) + \mathcal{O} (\gamma^2)
\Big\}
\end{align*}

Setting $\delta \!=\! \gamma$ and keeping the $\mathcal{O}(\gamma)$ terms yields: 
\begin{align*}
\dot\vz(t) + \frac{\gamma}{2}\ddot\vz(t) = - 3\alpha F(\vz(t)) + 3 \alpha \gamma J (\vz(t)) F(\vz(t)) \,,
\end{align*}

Re-writing the above in phase-space gives:
\begin{equation}\tag{LA3-GDA-HRDE}\label{eq:la3-gda_hrde3}
\begin{split}
  \dot{\vz}(t) & = \vomega(t) \\
  \dvom(t)     & = - \frac{2}{\gamma} \vomega(t) - \frac{6\alpha}{\gamma} F(\vz(t)) + 6 \alpha \cdot J(\vz(t)) \cdot F(\vz(t))  \,.
\end{split} 
\end{equation}

\paragraph{~\ref{eq:lookahead}4-GDA.}
For~\ref{eq:lookahead}4-GDA, replacing with \eqref{eq:predicted_gda_4} we get:
\begin{equation} 
\begin{aligned}\label{eq:la4_gda_interm1}
    \vz_{n+1} &= \vz_n + \alpha(\tilde \vz_{n+4} - \vz_n) \\
    &= \vz_n + \alpha\gamma\big[\underbrace{
    - F(\vz_n) - F\big(\vz_n - \gamma F(\vz_n)\big)
    - F\big( \vz_n - \gamma F(\vz_n) - \gamma F(\vz_n - \gamma F(\vz_n)) \big)}_{\text{same as for LA3-GDA}}\\
    &- \underbrace{F \Bigg(\vz_n
    - F(\vz_n) - F\big(\vz_n - \gamma F(\vz_n)\big)
    - F\big( \vz_n - \gamma F(\vz_n) - \gamma F(\vz_n - \gamma F(\vz_n)) \big)
    \Bigg)}_{(\star)}
    \big]
\end{aligned}
\end{equation}

Similarly as above for LA3-GDA, by performing consecutive TE in  parameter space for $(\star)$ we have:
\begin{align*}
(\star)= F(\vz_n) - 3\gamma J(\vz_n) F(\vz_n) + \mathcal{O} (\gamma^2)\,.
\end{align*}

Thus, using \ref{eq:time-taylor} as well as replacing the above in~\eqref{eq:la4_gda_interm1} we have: 
\begin{align*}
\frac{\dot{\vz}(n\delta)\delta \!+\! \frac{1}{2} \ddot{\vz}(n\delta)\delta^2  + \mathcal{O}(\delta^3)}{\gamma} 
&=
\alpha \Big\{
- 4 F(\vz(n\delta)) + 6\gamma J (\vz(n\delta)) F(\vz(n\delta)) + \mathcal{O} (\gamma^2)
\Big\}
\end{align*}

Setting $\delta \!=\! \gamma$ and keeping the $\mathcal{O}(\gamma)$ terms yields: 
\begin{align*}
\dot\vz(t) + \frac{\gamma}{2}\ddot\vz(t) = - 4\alpha F(\vz(t)) + 6 \alpha \gamma J (\vz(t)) F(\vz(t)) \,,
\end{align*}

Re-writing the above in phase-space gives:
\begin{equation}\tag{LA4-GDA-HRDE}\label{eq:la4-gda_hrde3}
\begin{split}
  \dot{\vz}(t) & = \vomega(t) \\
  \dvom(t)     & = - \frac{2}{\gamma} \vomega(t) - \frac{8\alpha}{\gamma} F(\vz(t)) + 12 \alpha \cdot J(\vz(t)) \cdot F(\vz(t))  \,.
\end{split} 
\end{equation}

\subsection{Summary: HRDEs of \ref{eq:lookahead}k-GDA}

Table~\ref{tab:la_gda_hrdes_summary} summarizes the obtained HRDEs for Lookahead-Minmax using GDA as a base optimizer, and generalizes it to any $k$.

\begin{table}[htb]
    \centering
    \begin{tabular}{r|l}
         GDA & 
         $\dot\vz(t) + \frac{\gamma}{2}\ddot\vz(t) = -  F(\vz(t))$  \\
         LA2-GDA & $\dot\vz(t) + \frac{\gamma}{2}\ddot\vz(t) = - 2\alpha \cdot F(\vz(t)) + \alpha \gamma \cdot J (\vz(t))  \cdot F(\vz(t))$ \\
         LA3-GDA & $\dot\vz(t) + \frac{\gamma}{2}\ddot\vz(t) = - 3\alpha  \cdot F(\vz(t)) + 3 \alpha \gamma  \cdot J (\vz(t))  \cdot F(\vz(t))$\\
         LA4-GDA & 
         $\dot\vz(t) + \frac{\gamma}{2}\ddot\vz(t) = - 4\alpha  \cdot F(\vz(t)) + 6 \alpha \gamma  \cdot J (\vz(t)) \cdot  F(\vz(t))$\\
         LA5-GDA & 
         $\dot\vz(t) + \frac{\gamma}{2}\ddot\vz(t) = - 5\alpha  \cdot F(\vz(t)) + 10 \alpha \gamma  \cdot J (\vz(t)) \cdot  F(\vz(t))$ \\
         \multicolumn{1}{c}{$\vdots$} & \multicolumn{1}{c}{$\vdots$} \\
         LA$k$-GDA & 
         $\dot\vz(t) + \frac{\gamma}{2}\ddot\vz(t) = - k \alpha  \cdot F(\vz(t)) + (\sum_{i=1}^{k-1} i)\cdot \alpha \gamma  \cdot J (\vz(t)) \cdot  F(\vz(t))$\\[.5em]
    \end{tabular}
    \caption{Summary: HRDEs of LAk-GD.}
    \label{tab:la_gda_hrdes_summary}
\end{table}
\section{Analysis of Lookahead Dynamics}
We first prove general convergence of the Lookahead dynamics for the \ref{eq:bilinear_game} in the discrete domain. We write the phase-space representation of the differential equation, then enforce convergence on the eigen values of the coefficient matrix using the Routh-Hurwitz criterion, as introduced in ~\S \ref{sec:prelim}.
\label{sec:theorem1-proof}
\textbf{Proof of Theorem \ref{theo:convergence}}

For the LA - kGD optimiser, we have the following differential equation
\begin{align*}
    \dot{x}(t) = \omega(t)
\end{align*}
\begin{align*}    
    \dot{\omega}(t) = - \frac{2}{\gamma} \omega(t) - \frac{2k}{\gamma} F(z(t)) + \frac{k(k -1)}{\gamma } \alpha \gamma J(z(t)) \cdot F(z(t))
\end{align*}

By denoting $\dot{x}(t) = \omega_x(t)$ and $\dot{y}(t) = \omega_y(t)$ we get the following

\begin{align*}
    \begin{bmatrix}
        \dot{x}(t) \\ 
        \dot{y}(t) \\ 
        \dot{\omega}_x(t) \\ 
        \dot{\omega}_y(t) 
    \end{bmatrix} = \underbrace{\begin{bmatrix}
        0 && 0 && \mathbf{I} && 0 \\
        0 && 0 &&  0 && \mathbf{I} \\
        - k (k -1) \alpha \vA^2 && -\dfrac{2k \alpha \vA}{ \gamma} && - \frac{2}{\gamma} \mathbf{I} && 0 \\
        \dfrac{2k \alpha \vA}{\gamma} && - k(k -1) \alpha \vA^2 && 0 && -\frac{2}{\gamma} \mathbf{I} 
    \end{bmatrix}}_{\triangleq \hspace{3pt} \mathcal{C}_{LA-kGD}} \cdot \\ \begin{bmatrix}
        x(t) \\
        y(t) \\
        \omega_x(t) \\
        \omega_y(t)
    \end{bmatrix}
\end{align*}
To obtain the eigenvalues $\lambda \in \mathbf{C}$ of $\mathbf{C}_{LA-kGD}$, we have

$\det(\mathcal{C}_{LA-kGD} - \lambda\mathbf{I})$
\begin{align*}
      &= \det\Bigg( 
    \begin{bmatrix}
         -\lambda\mathbf{I} & 0 & \mathbf{I} & 0 \\
        0 & -\lambda \mathbf{I} &  0 & \mathbf{I} \\
        - k (k -1) \alpha \vA^2 & -\dfrac{2k \alpha \vA}{ \gamma} & - (\frac{2}{\gamma} +\lambda) \mathbf{I} & 0 \\
        \dfrac{2k \alpha \vA}{ \gamma} & - k(k -1) \alpha \vA^2 & 0 & -(\frac{2}{\gamma } + \lambda)  \mathbf{I}
    \end{bmatrix}  
    \Bigg)  \\
    &=  \det \Bigg( 
    \begin{bmatrix}
        \lambda(\frac{2}{\gamma} + \lambda) & 0\\
        0 & \lambda(\frac{2}{\gamma} + \lambda)
    \end{bmatrix} 
    - \underbrace{\begin{bmatrix}
        - k (k -1) \alpha \vA^2& -\dfrac{2k \alpha \vA}{\gamma} \\ 
         \dfrac{2k \alpha \vA}{ \gamma} & - k(k -1) \alpha \vA^2
    \end{bmatrix}}_{\triangleq \mathbf{D}} 
    \Bigg) .
\end{align*}

Let $\mu = \mu_1 + i \mu_2  \in \mathbb{C}$ denote the eigenvalues of $\mathbf{D}$. The characteristic equation becomes
\[
\lambda\left( \frac{2}{\gamma} + \lambda \right) - \mu = 0.
\]
Let $\beta = \frac{2}{\gamma}$. Using the generalized Hurwitz theorem for polynomials with complex coefficients, we construct the following generalized Hurwitz array:

\[
\begin{array}{c|ccc}
\lambda^2 & 1 & 0 & \mu_1 \\
\lambda^1 & \beta & \mu_2 & 0 \\
 & -\mu_2 & \beta \mu_2 & 0 \\
   \lambda^0      & -\mu_2^2 - \beta^2 \mu_1 & 0 & 0 \\
\end{array}
\]

The sign of the last row determines the stability of the polynomial. Since $\beta > 0$, the system is stable if and only if:
\begin{equation}
    \label{eq:stability-equation}
    \mu_1 < -\frac{1}{\beta^2} \mu_2^2 \,.
\end{equation}

Thus, it suffices to show:
\begin{equation}
\label{eq:stability-condition}
\Re(\mu(z)) < -\frac{1}{\beta^2} \left( \Im(\mu(z)) \right)^2 \,.
\end{equation}

We have that: 
\begin{align*}
    \mu(z) = \bar{z}^T D z &= [\bar{x}^T \quad \bar{y}^T]\begin{bmatrix}
        - k (k -1) \alpha \vA^2 && -\dfrac{2k \alpha \vA}{ \gamma} \\ 
         \dfrac{2k \alpha \vA}{ \gamma} && - k(k -1) \alpha \vA^2
    \end{bmatrix}\begin{bmatrix}
        x \\
        y
    \end{bmatrix} \\
    &= -k(k -1)\alpha ( \|\vA x\|_2^2+ \|\vA y\|_2^2) + k \alpha (\frac{\bar{y}^T \vA x}{ \gamma} - \frac{\bar{x}^T\vA y}{\gamma}) \,.
\end{align*}

Note that $\bar{x}^T\vA y$ is the complex conjugate of $\bar{y}^T \vA x$, hence $\bar{y}^T \vA x - \bar{x}^T\vA y = 2i \Im(\bar{x}^T \vA y) $

Using this in \ref{eq:stability-equation}, we get
\begin{align*}
    \|\vA x\|_2^2 + \|\vA y\|_2^2 \geq 4 \dfrac{\alpha k}{ (k -1)} \|\bar{x}^T \vA y\|^2   \,.
\end{align*}

\section{Convergence Analysis in Continuous Time}
In this section, we extend the analysis of the discrete-time algorithms in the continuous time. We take the Laplace transform (\ref{eq:Laplace}) of the HRDE, and reuse the Routh-Hurwitz criterion to enforce convergence on the poles of the frequency duals. We prove the divergence of \ref{eq:gda-hrde} and convergence of \ref{eq:la-hrde} on \ref{eq:bilinear_game} as discussed in ~\S \ref{sec:convergence}.  

\subsection{Showing Divergence of \ref{eq:gda-hrde} for \ref{eq:bilinear_game} Through \eqref{eq:Laplace}}\label{app:gda_convergence_analysis}

 As shown in~\citep{chavdarova2023hrdes}, Gradient decent's $\mathcal{O}(\gamma)$-HRDE~\citep{chavdarova2023hrdes} is:
\begin{equation}
\tag{GD-HRDE}
    \ddot{\vz}(t)     = - \frac{2}{\gamma}\cdot \dot{\vz}(t) - \frac{2}{\gamma} \cdot F(\vz(t)) \,,
\label{eq:gda-hrde}
\end{equation}
where $\vz$ represents the vector of players $\vz(t) \triangleq (\vx(t) , \vy(t))^\intercal\,.$
$\vz(t) = \begin{bmatrix}
    \vx(t) \\
    \vy(t)
\end{bmatrix}$ 
and $F(\cdot)$ is the operator defined in \eqref{eq:vi}.

We compute the gradient field and the jacobian of the parameterized bilinear game \ref{eq:bilinear_game}. Note that we are using the continuous time representation of the two players.
\begin{equation}
    z(t) = \begin{bmatrix}
        \vx(t) \\
        \vy(t)
    \end{bmatrix}
\end{equation}
\begin{equation}
    F(\vz(t)) = \begin{bmatrix}
     \vA\vy(t) \\
     - \vA\vx(t)\end{bmatrix}
\end{equation}
We rewrite the  \ref{eq:gda-hrde} for the 2 players as follows
\begin{equation}
\label{eq:hrde-gda-x}
    \ddot{\vx}(t) =  - \frac{2}{\gamma} \cdot \vA\dot{\vx}(t) - \frac{2}{\gamma} \cdot \vy (t)  
\end{equation}
\begin{equation}
\label{eq:hrde-gda-y}
    \ddot{\vy}(t) = - \frac{2}{\gamma} \cdot \vA \dot{\vy}(t) + \frac{2}{\gamma} \cdot \vx(t)  
\end{equation}
We take the Laplace transform of the \ref{eq:hrde-gda-x} and \ref{eq:hrde-gda-y}

\begin{equation}
       \mathbf{X}(s) =  -\dfrac{2 }{ \gamma (s^2 I+ \frac{2\vA s}{\gamma})} \mathbf{Y}(s) + \dfrac{\dot{\vx}(0)+ \frac{\vA}{\gamma} \vx(0)}{s^2 I + \frac{2\vA s}{\gamma} } + \dfrac{(s + \frac{A}{\gamma})\vx(0)}{ s^2 + \frac{2\vA s}{\gamma} }
        \label{eq:freq_gda_X}
\end{equation}
\begin{equation}
         \mathbf{Y}(s) = -\dfrac{2}{ \gamma (s^2 I+ \frac{2\vA s}{\gamma})} \mathbf{X}(s) + \dfrac{\dot{\vy}(0)+ \frac{\vA}{\gamma} \vy(0)}{s^2 + \frac{2\vA s}{\gamma} } + \dfrac{(s + \frac{\vA}{\gamma})\vy(0)}{ s^2 + \frac{2\vA s}{\gamma} }
        \label{eq:freq__gda_Y}
\end{equation}
Here, $s \in \mathbb{C}$ represents the complex frequency, $\mathbf{X}(s) =  \int_{0}^{\infty} x(t) e^{-st} \, dt $, $\mathbf{Y}(s) =  \int_{0}^{\infty} \vy(t) e^{-st} \, dt$ are the transfer functions of the trajectories $\vx(t)$ and $\vy(t)$. We assume that $\vx(t)$ and $\vy(t)$ are defined for $t>0$ i.e. a unilateral
Laplace transform. .

 Since $\vX(s)$ and $\vY(s)$ are linear equations, we compute the combined transform function for $\vX(s)$. The poles of the combined transform function are $s \in \mathbf{C}$ such that
 \begin{equation}
 \label{eq:gd_condition_roots}
    det(s^4I + \frac{2s^3}{\gamma} (I+\vA) + \frac{4s^2}{\gamma^2} \vA - \frac{4}{\gamma^2}I) = 0
 \end{equation}

 We use Routh Hurwitz criterion to analyse the stability of \ref{eq:gd_condition_roots}. Let $\lambda_i \in \Lambda$ be eigenvalue of A.
 
 Since the determinant is zero when any eigenvalue of the matrix polynomial is zero, consider an eigenvalue \( \lambda_i \) of \(\mathbf{A}\). The scalar polynomial corresponding to \(\lambda_i\) is:

\begin{equation}
p_i(s) = s^4 + \frac{2 s^3}{\gamma} (1 + \lambda_i) + \frac{4 s^2}{\gamma^2} \lambda_i - \frac{4}{\gamma^2} = 0.
\end{equation}

We rewrite \(p_i(s)\) as

\begin{equation}
\label{eq:general-polynomial}
p_i(s) = s^4 + a_3 s^3 + a_2 s^2 + a_1 s + a_0,
\end{equation}

where

\begin{equation}
a_3 = \frac{2}{\gamma} (1 + \lambda_i), \quad
a_2 = \frac{4}{\gamma^2} \lambda_i, \quad
a_1 = 0, \quad
a_0 = -\frac{4}{\gamma^2}.
\end{equation}

The Routh array for \ref{eq:general-polynomial} is:

\[
\begin{array}{c|cc}
s^4 & 1 & a_2 \\
s^3 & a_3 & a_1 \\
s^2 & b_1 & a_0 \\
s^1 & c_1 & 0 \\
s^0 & a_0 & 0
\end{array}
\]

where

\begin{equation}
b_1 = \frac{a_3 a_2 - 1 \cdot a_1}{a_3} = \frac{a_3 a_2}{a_3} = a_2,
\end{equation}

(since \(a_1=0\)),

\begin{equation}
c_1 = \frac{b_1 a_1 - a_3 a_0}{b_1} = \frac{0 - a_3 a_0}{a_2} = -\frac{a_3 a_0}{a_2}.
\end{equation}

The Routh-Hurwitz stability conditions are as follows

\begin{equation}
1 > 0, \quad
a_3 = \frac{2}{\gamma} (1 + \lambda_i) > 0, \quad
b_1 = a_2 = \frac{4}{\gamma^2} \lambda_i > 0,
\end{equation}

\begin{equation}
c_1 = -\frac{a_3 a_0}{a_2} = -\frac{\frac{2}{\gamma} (1 + \lambda_i) \cdot \left(-\frac{4}{\gamma^2}\right)}{\frac{4}{\gamma^2} \lambda_i} = \frac{2 (1 + \lambda_i)}{\gamma} \cdot \frac{1}{\lambda_i} > 0,
\end{equation}

and finally

\begin{equation}
a_0 = -\frac{4}{\gamma^2} > 0.
\end{equation}

\begin{equation}
a_0 = -\frac{4}{\gamma^2} > 0 \implies -\frac{4}{\gamma^2} > 0 \implies \frac{1}{\gamma^2} < 0,
\end{equation}

which is false. Hence the system is unstable, and GDA diverges for \ref{eq:bilinear_game}, which confirms the expected behavior of GDA for games.

\subsection{Convergence of \ref{eq:la-hrde} for 
\ref{eq:bilinear_game} using \eqref{eq:Laplace}}
\label{app:la_bg_convergence}
We first present the proof of Theorem \ref{theo:2} by taking the inverse Laplace Transform (\ref{inverse-laplace}) of the frequency dual of \ref{eq:la-hrde}. We then perform a similar analysis on the poles of the frequency dual as described in \ref{app:gda_convergence_analysis} to find the convergence criteria for Lookahead.

\textbf{Proof of Theorem \ref{theo:2}}
We compute similarly, the gradient field and the Jacobian of the bilinear game \eqref{eq:bilinear_game}. 
\begin{equation}
    \vz(t) = \begin{bmatrix}
        \vx(t) \\
        \vy(t)
    \end{bmatrix}
\end{equation}

\begin{equation}
    F(\vz(t)) = \begin{bmatrix}
    \vA \vy(t) \\
    -\vA \vx(t)\end{bmatrix}
\end{equation}
\begin{equation}
\label{eq:operator}
    \mathcal{J} F(\vz(t))  = \begin{bmatrix}
        0 && \vA  \\
        -\vA && 0
    \end{bmatrix}
\end{equation}
\begin{equation}
\label{eq:product}
    \mathcal{J} F(\vz(t)) \cdot F(\vz(t)) = \begin{bmatrix}
       -\vA^2 \vy(t)\\
         -\vA^2 \vx(t)
    \end{bmatrix}
\end{equation}
We rewrite the  \ref{eq:la-hrde} for the 2 players using the equations \ref{eq:operator} and \ref{eq:product}

\begin{equation}
\begin{split}
\label{eq:substituted-hrde-x}
        \ddot{\vx}(t) &= -\frac{2}{\gamma} \dot{\vx}(t) - \frac{2k\alpha}{\gamma} \vA \vy(t) - \alpha k(k-1) \vA^2 \vx(t)  
\end{split}
\end{equation}
\begin{equation}
\begin{split}
\label{eq:substituted-hrde-y}
    \ddot{\vy}(t) &= -\frac{2}{\gamma} \dot{\vy}(t) + \frac{2k\alpha}{\gamma} \vA \vx(t) - 2\alpha k(k-1) \vA^2 \vy(t)
\end{split}
\end{equation}

We take the Laplace transform of \ref{eq:substituted-hrde-x} and \ref{eq:substituted-hrde-y}.
\begin{equation}
       \left(s^2 I + \frac{2}{\gamma} s I + \alpha k(k-1)  \vA^2\right) \vX(s) + \frac{2k\alpha}{\gamma} \vA \vY(s) = \frac{2}{\gamma} \vx(0) + s\vx(0) + \dot{\vx}(0)
        \label{eq:freq_X}
\end{equation}
\begin{equation}
         -\frac{2k\alpha}{\gamma} \vA \vX(s) + \left(s^2 I + \frac{2}{\gamma} s I + \alpha k(k-1)  \vA^2\right) \vY(s) = \frac{2}{\gamma} \vy(0) + s\vy(0) + \dot{\vy}(0)
        \label{eq:freq_Y}
\end{equation}

Taking the Inverse Laplace transform (\ref{inverse-laplace}) of the above yields the solution equations in \ref{eq:trajectory_1}.
We solve the system \eqref{eq:freq_X}--\eqref{eq:freq_Y} by applying the eigendecomposition $\vA = \mU \mLambda \mU^\intercal$ with $\mLambda = \mathrm{diag}(\lambda_i)$ and orthogonal $\mU$. Left-multiplying by $\mU^\intercal$ and defining $\vx(s) = \mU \vx_{\text{eig}}(s)$ and $\vy(s) = \mU \vy_{\text{eig}}(s)$, the system decouples along each eigendirection $i$:

\begin{align}
\left(s^2 + \frac{2}{\gamma} s + \alpha k(k{-}1) \lambda_i^2 \right) x_i(s) + \frac{2k\alpha}{\gamma} \lambda_i y_i(s) &= s x_i(0) + \dot{x}_i(0) + \frac{2}{\gamma} x_i(0), \\
-\frac{2k\alpha}{\gamma} \lambda_i x_i(s) + \left(s^2 + \frac{2}{\gamma} s + \alpha k(k{-}1) \lambda_i^2 \right) y_i(s) &= s y_i(0) + \dot{y}_i(0) + \frac{2}{\gamma} y_i(0).
\end{align}

Solving this $2 \times 2$ system gives expressions of the form:
\[
x_i(s),\ y_i(s) = \frac{\text{linear in } s}{(s + \frac{1}{\gamma})^2 - \omega_i^2}, \quad \text{where } \omega_i^2 = \alpha k(k{-}1)\lambda_i^2 - \frac{1}{\gamma^2}.
\]

Applying the standard Laplace inverses:
\begin{align*}
\mathcal{L}^{-1} \left\{ \frac{1}{(s + \frac{1}{\gamma})^2 - \omega_i^2} \right\} &= e^{-t/\gamma} \frac{\sinh(\omega_i t)}{\omega_i}, \\
\mathcal{L}^{-1} \left\{ \frac{s + \frac{1}{\gamma}}{(s + \frac{1}{\gamma})^2 - \omega_i^2} \right\} &= e^{-t/\gamma} \cosh(\omega_i t),
\end{align*}

we obtain the time-domain solutions:
\[
x_i(t) = e^{-t/\gamma} \left[ \cosh(\omega_i t)\, x_i(0) + \frac{1}{\omega_i} \sinh(\omega_i t) \left( \dot{x}_i(0) + \frac{1}{\gamma} x_i(0) \right) \right] - \frac{2k\alpha}{\gamma} (g_i * y_i)(t),
\]
\[
y_i(t) = e^{-t/\gamma} \left[ \cosh(\omega_i t)\, y_i(0) + \frac{1}{\omega_i} \sinh(\omega_i t) \left( \dot{y}_i(0) + \frac{1}{\gamma} y_i(0) \right) \right] + \frac{2k\alpha}{\gamma} (g_i * x_i)(t),
\]

where $g_i(t) = e^{-t/\gamma} \frac{\sinh(\omega_i t)}{\omega_i}$, $\mG(t) = \mU\, \mathrm{diag}(g_i(t))\, \mU^\intercal$, and $*$ denotes convolution. Finally, transforming back to the original coordinates gives the full vector solution:

\begin{align}
\vx(t) &= -\frac{2k\alpha}{\gamma} (\mG * \vy)(t) + \mU\, \mathrm{diag}\left( e^{-t/\gamma} \left[ \cosh(\omega_i t)\, x_i(0) + \frac{1}{\omega_i} \sinh(\omega_i t) \left( \dot{x}_i(0) + \frac{1}{\gamma} x_i(0) \right) \right] \right) \mU^\intercal, \\
\vy(t) &= \frac{2k\alpha}{\gamma} (\mG * \vx)(t) + \mU\, \mathrm{diag}\left( e^{-t/\gamma} \left[ \cosh(\omega_i t)\, y_i(0) + \frac{1}{\omega_i} \sinh(\omega_i t) \left( \dot{y}_i(0) + \frac{1}{\gamma} y_i(0) \right) \right] \right) \mU^\intercal,
\end{align}

\textbf{Note}: For a purely potential game, the terms of $\vY(s)$ cancel out and simplifies the general \ref{eq:trajectory_1} to:

\begin{equation*}
     \vx(t) = \hspace{2pt}  \mU \, \mathrm{diag}\left( 
  e^{-\frac{t}{\gamma}} \left[ 
  \cosh(\omega_i t) \, \vx_i(0) + \frac{1}{\omega_i} \sinh(\omega_i t) \left( \dot{\vx}_i(0) + \frac{1}{\gamma} \vx_i(0)
  \right)
  \right]\right) \mU^\intercal  \,.
\end{equation*}

which confirms that in pure minimization settings, the solution $\vx(t)$ does not depend on $\vy(t)$.

\textbf{Convergence Analysis}
We substitute \ref{eq:freq_Y} in \ref{eq:freq_X} to get the joint Laplace transform for $\vx(t)$ as follows 
\begin{equation*}
\label{eq:freq_hrde}
\tag{$X$--TF}
\begin{aligned}
\mathbf{X}(s) =\; & \dfrac{-2k\gamma\alpha \vy(0)\, (s+ \frac{2}{\gamma})\vA - 2k\alpha\gamma\, \dot{\vy}(0)\vA}{\gamma^2 \left(s^2 + \frac{2s}{\gamma} + \alpha k(k{-}1) \vA^2\right)^2 + 4k^2\alpha^2} \\
& + \dfrac{\dot{\vx}(0) + \left(s+\frac{2}{\gamma}\right)\vx(0)}{\gamma^2 \left(s^2 + \frac{2s}{\gamma} + \alpha k(k{-}1) \vA^2\right)^2 + 4k^2\alpha^2} \cdot \gamma^2\left(s^2 + \frac{2s}{\gamma} + \alpha k(k{-}1)\right)
\end{aligned}
\end{equation*}
\begin{equation*} \label{eq:freq_hrde-y} \tag{$Y$--TF}
\begin{aligned}
    \mathbf{Y}(s) =\; & \dfrac{2k\gamma\alpha \vx(0) (s+ \frac{2}{\gamma}) \vA +2k\alpha\gamma\dot{\vx}(0) \vA}{\gamma^2 (s^2 +\frac{2s}{\gamma}+\alpha k(k-1) \vA^2)^2 + 4k^2\alpha^2} +\\
& \dfrac{\dot{\vy}(0) + (s+\frac{2}{\gamma})\vy(0)}{\gamma^2 (s^2 +\frac{2s}{\gamma}+ \alpha k(k-1) \vA^2 )^2 + 4k^2\alpha^2} \gamma^2(s^2 +\frac{2s}{\gamma}+\alpha k(k-1))
\end{aligned}
\end{equation*}

The characteristic equation is then: 
\begin{equation}
    \left( s^2 + \frac{2}{\gamma}s + \alpha k(k-1) \vA^2 \right)^2 + \left( \frac{2k\alpha}{\gamma} \vA \right)^2 = 0 \,.
\end{equation}
We use the Routh-Hurwitz criterion to analyse the coefficients of the characteristic equations. We arrive at the following convergence condition (using similar steps as in \ref{app:gda_convergence_analysis})

\begin{equation}
    \label{eq:convergence-condition}
    \tag{BG-Cond}
    \alpha < \frac{k-1}{k} \,.
\end{equation}
\subsection{Convergence of \ref{eq:la-hrde} for \ref{eq:quadratic_game} using \eqref{eq:Laplace}}
\label{sec:convergence-la-qd}
In this section, we use our framework to analyse the convergence of Lookahead on \ref{eq:quadratic_game}. We use the same method of using \ref{eq:Laplace} to get the frequency duals, then enforce convergence using the roots of the characteristic equation. 

The operator for the \ref{eq:quadratic_game} game $F: \mathbb{R}^n \rightarrow \mathbb{R}^n$, is as follows,

\begin{align*}
    F = \begin{bmatrix}
        \vC\vx + \vA\vy \\
        -\vA\vx + \vB\vy
    \end{bmatrix}
\end{align*}

We rewrite $F$ as the following 
\begin{align*}
    F = \begin{bmatrix}
        \vC\\ 
        -\vA
    \end{bmatrix} x + \begin{bmatrix}
        \vA \\
        \vB
    \end{bmatrix} y
\end{align*}

\begin{align*}
    F = \vT_x \vx + \vT_y \vy
\end{align*}

where $\vT_x = \begin{bmatrix}
    \vC \\ 
    -\vA
\end{bmatrix}$ and $\vT_y = \begin{bmatrix}
    \vA \\ 
    \vB
\end{bmatrix}$ are transformation matrices for notational convenience.

Additionally, 

\begin{align*}
    \mathcal{J}F = \begin{bmatrix}
        \vC && \vA \\ 
        -\vA && \vB
    \end{bmatrix}
\end{align*}

\begin{align*}
    \mathcal{J}F \cdot F &= \begin{bmatrix}
        \vC & \vA \\ 
        -\vA & \vB
    \end{bmatrix} \begin{bmatrix}
        \vC\vx + \vA\vy \\
        -\vA\vx + \vB\vy
    \end{bmatrix}\\ 
    &= \begin{bmatrix}
        (\langle \vC, \vC \rangle - \langle \vA, \vA \rangle) \vx + (\langle \vC, \vA \rangle + \langle \vA, \vB \rangle) \vy \\ 
        - (\langle \vA, \vC \rangle + \langle \vB, \vA \rangle) \vx + (\langle \vB, \vB \rangle - \langle \vA, \vA \rangle) \vy 
    \end{bmatrix}\\ 
    &= \begin{bmatrix}
        \langle \vC, \vC \rangle - \langle \vA, \vA \rangle  \\ 
        - \langle \vA, \vC \rangle - \langle \vB, \vA \rangle 
    \end{bmatrix} \cdot \vx + \begin{bmatrix}
        \langle \vC, \vA \rangle - \langle \vA, \vB \rangle  \\ 
        \langle \vB, \vB \rangle - \langle \vA, \vA \rangle 
    \end{bmatrix} \cdot \vy \\
    &= \vT^{'}_x \vx + \vT^{'}_y \vy
\end{align*}

Hence, the HRDE is as follows

\begin{align*}
    \dot{\vx}(t) + \frac{1}{2} \gamma \ddot{\vx}(t) &= -k \alpha (\vC \vx(t) + \vA \vy(t)) \\ 
    &\quad+ k(k - 1)\alpha \gamma  \big((\langle \vC, \vC \rangle - \langle \vA, \vA \rangle) \vx(t) +  (\langle \vC, \vA \rangle - \langle \vA, \vB \rangle) \vy(t)\big),
\end{align*}

\begin{align*}
    \dot{\vy}(t) + \frac{1}{2}\gamma\ddot{\vy}(t) = & -k \alpha (-\vA\vx(t) + \vB\vy(t)) \\
    & + k(k - 1)\alpha \gamma \left( - (\langle \vA, \vC \rangle + \langle \vB, \vA \rangle) \vx(t) + (\langle \vB, \vB \rangle - \langle \vA, \vA \rangle) \vy(t) \right)
\end{align*}

We introduce the real valued constants $m_x, m_y, l_x, l_y$ for convience. We then take the Laplace transform of the above differential equations, perform simple linear substitution to obtain the characteristic polynomial of the form 

\begin{align}
    (\frac{1}{2} \gamma s^2 + s + m_x) (\frac{1}{2} \gamma s^2 + s + m_y) + l_{x} l_{y} = 0   
\end{align}

where

\begin{align*}
    m_x = k \alpha \cdot [1 \quad 0] \vT_x - k(k - 1)\alpha \gamma \cdot [1 \quad 0] \vT_x^{'}
\end{align*}

\begin{align*}
    m_y = k \alpha \cdot [0 \quad 1] \vT_y - k(k - 1)\alpha \gamma \cdot [0 \quad 1] \vT_y^{'}
\end{align*}

\begin{align*}
    l_x = -k \alpha \cdot [1 \quad 0] \vT_y + k(k - 1)\alpha \gamma \cdot [1 \quad 0] \vT_y^{'}
\end{align*}

\begin{align*}
    l_y = -k \alpha \cdot [0 \quad 1] \vT_x + k(k - 1)\alpha \gamma \cdot [0 \quad 1] \vT_x^{'}
\end{align*}

Following the Routh - Hurowitz criterion, the condition for convergence is as follows 

\begin{align*}
    \gamma (m_x - m_y)^2 + 4 (m_x + m_y) - 4l_xl_y \gamma > 0
\end{align*}
 The full convergence condition is as follows
 \begin{equation}
\label{eq:convergence-inequality-full}\tag{QD-Cond}
\begin{split}
& \left(
\gamma \cdot \left[ \begin{matrix} 1 & 0 \end{matrix} \right] \left( k \alpha T_x - k(k - 1)\alpha \gamma T_x' \right) \right.   \left. - \gamma \cdot \left[ \begin{matrix} 0 & 1 \end{matrix} \right] \left( k \alpha T_y - k(k - 1)\alpha \gamma T_y' \right) \right)^2 \\
& +\, 4\gamma \cdot \left( 
\left[ \begin{matrix} 1 & 0 \end{matrix} \right] \left( k \alpha T_x -k(k - 1)\alpha \gamma T_x' \right) \right.  \left. + \left[ \begin{matrix} 0 & 1 \end{matrix} \right] \left( k \alpha T_y - k(k - 1)\alpha \gamma T_y' \right) \right) \\
& -\, 4\gamma^2 \cdot \left( 
\left[ \begin{matrix} 1 & 0 \end{matrix} \right] \left( -k \alpha T_y + k(k - 1)\alpha \gamma T_y' \right) \right.   \left. \cdot \left[ \begin{matrix} 0 & 1 \end{matrix} \right] \left( -k\alpha T_x + k(k - 1)\alpha \gamma T_x' \right) \right) \\
& \hspace{200pt} \quad > 0
\end{split}
\end{equation}
For a purely Hamiltonian game, $\vB = \vC = \mathbf{0} _{\frac{d}{2} \times \frac{d}{2}}$; that is~\ref{eq:bilinear_game}, and the convergence condition \ref{eq:convergence-inequality-full} reduces to \ref{eq:cond_hamiltonian}.

\section{\texorpdfstring{$O(\gamma^2)$ High Resolution Differential Equation for Lookahead}{O(gamma²) High Resolution Differential Equation for Lookahead}}\label{app: hrde-derive-second}

In this section, we derive a second order differential equation for Lookahead, as described in ~\S \ref{sec:ogamma-2}. 
We follow a similar method as that of the $O(\gamma)$ HRDE in Appendix \ref{app:la_hrde_derivation}, except also include second order terms in the differential equation, and take the limit as $O(\gamma^3) \rightarrow 0$.

\textbf{LA2-GDA}
The iterates of LA2-GDA are obtained as follows:
\begin{align*}
    \vz_{n+1} = \vz_n + \alpha (\tilde\vz_{n+2} - \tilde\vz_n) = \vz_n + \alpha \Big(
    -\gamma F(\vz_n) - \gamma F\big(
    \vz_n - \gamma F(\vz_n)
    \big)
    \Big)\,.
\end{align*}

Using TE (including terms upto $O(\gamma^2)$ for the right hand side, we get the following
\begin{align*}
\label{LA2-TE}
\tag{LA2-TE}
    F(\vz_n  - \gamma F(\vz_n)) = F(\vz_n) - \gamma J(\vz_n)F(\vz_n) + \frac{\gamma^2}{2} F(\vz_n)^T H(\vz_n) F(\vz_n) + \mathcal{O}(\gamma^3)
\end{align*}

Using TE, we get (where $\delta$ and $\gamma$ are the step sizes in time and parameter space, resp.):

\begin{align*}
    \frac{\dot \vz (n\delta) + \frac{1}{2} \delta^2 \ddot{\vz} (n\delta) + \frac{1}{6} \delta^3 \dddot{z}(n\delta) + \mathcal{O} (\delta^4) }{\gamma}
    = \alpha \Big( - 2F(\vz (n\delta)) + \gamma J(\vz(t))F(\vz(n\delta)) - \\  \frac{\gamma^2}{2} F(\vz_n)^T H(\vz_n) F(\vz_n)
    \Big) \,.
\end{align*}
Setting $\delta \!=\! \gamma$ and keeping the $\mathcal{O}(\gamma^2)$ terms yields:
\begin{align*}
\tag{LA2-HRDE}
\label{LA2-HRDE}
\dot \vz (t) + \frac{1}{2} \gamma \ddot{\vz} (t) + \frac{1}{6} \gamma^2 \dddot{z}(t) 
    = \alpha \Big( - 2F(\vz (t)) + \gamma J(\vz(t))F(\vz(t)) - \\ \frac{\gamma^2}{2} F(\vz(t))^T H(\vz(t)) F(\vz(t))
    \Big) \,.
\end{align*}

\textbf{LA3-GDA}
For LA3-GDA we have:
\begin{equation}
\begin{aligned}
    \vz_{n+1} &= \vz_n + \alpha(\tilde \vz_{n+3} - \vz_n) \\
    &= \vz_n + \alpha\gamma\big[
    - F(\vz_n) - F\big(\vz_n - \gamma F(\vz_n)\big)
    \\& - \underbrace{F\big( \vz_n - \gamma F(\vz_n) - \gamma F(\vz_n - \gamma F(\vz_n)) \big)}_{(\star)}
    \big]
\end{aligned}
\end{equation}
Similarly, by doing TE in coordinate space for ($\star$),  we get:
\begin{align*}
    &F\big( \vz_n - \gamma F(\vz_n) - \gamma F(\vz_n - \gamma F(\vz_n)) \big) \\
    &= F(\vz_n) + J(\vz_n)\big(-\gamma F(\vz_n) - \gamma F(\vz_n - \gamma F(\vz_n))\big) \\
    &\quad + \frac{1}{2} \big(-\gamma F(\vz_n) - \gamma F(\vz_n - \gamma F(\vz_n))\big)^T \\
    &\quad\quad H(\vz_n) \big(-\gamma F(\vz_n) - \gamma F(\vz_n - \gamma F(\vz_n))\big)
\end{align*}

Using \ref{LA2-TE}, and keeping only $\mathcal{O}(\gamma^2)$ terms, we get

\begin{align*}
    F\big( \vz_n - \gamma F(\vz_n) - \gamma F(\vz_n - \gamma F(\vz_n)) \big)& = F(\vz_n) \\& +  J(\vz_n)( -2\gamma F(\vz_n) + \gamma^2 J(\vz_n)F(\vz_n)) \\&+ \frac{1}{2}(2\gamma F(\vz_n))^T H(\vz_n) (2\gamma F(\vz_n)) \\
\end{align*}

\begin{align*}
    \label{LA3-TE}
    F\big( \vz_n - \gamma F(\vz_n) - \gamma F(\vz_n - \gamma F(\vz_n)) \big) &= F(\vz_n) - 2\gamma J(\vz_n)F(\vz_n) \\& + \gamma^2 (J(\vz_n))^2 F(\vz_n) + \\& 2\gamma^2 F(\vz_n)^T H(\vz_n) F(\vz_n)
    \tag{LA3-TE}
\end{align*}

\begin{equation}     
    \begin{aligned}         
        \vz_{n+1}  - \vz_n  &= + \alpha \gamma \big[ -F(\vz_n) - F(\vz_n) + \gamma J(\vz_n)F(\vz_n) 
        - \frac{\gamma^2}{2} F(\vz_n)^T H(\vz_n) F(\vz_n) \\
        &\quad - F(\vz_n) + 2\gamma J(\vz_n)F(\vz_n) 
        - \gamma^2 (J(\vz_n))^2 F(\vz_n) 
        \\& - 2\gamma^2 F(\vz_n)^T H(\vz_n) F(\vz_n) \big]     
    \end{aligned} 
\end{equation}

Using TE for the left hand side, we get
\begin{equation}
    \begin{aligned}
        \frac{\dot \vz (t) + \frac{1}{2} \delta^2 \ddot{\vz} (t) + \frac{1}{6} \delta^3 \dddot{\vz}(t) + \mathcal{O}(\delta^4)}{\gamma}
        &= \alpha \Big( -3F(\vz_n) + 3\gamma J(\vz_n)F(\vz_n) \\
        &\quad - \frac{5}{2} \gamma^2 F(\vz_n)^T H(\vz_n)F(\vz_n) - \gamma^2 (J(\vz_n))^2 F(\vz_n) \Big)
    \end{aligned}
\end{equation}
Setting $\delta \!=\! \gamma$ and keeping the $\mathcal{O}(\gamma^2)$ terms yields:
\begin{align*}
\dot \vz (t) + \frac{1}{2} \gamma \ddot{\vz} (t) + \frac{1}{6} \gamma^2 \dddot{\vz}(t) 
&= \alpha \Big( -3F(\vz(t)) + 3\gamma J(\vz(t))F(\vz(t)) \\
&\quad - \frac{5}{2} \gamma^2 F(\vz(t))^T H(\vz(t))F(\vz(t)) 
       \\ & \quad - \gamma^2 J(\vz(t))^2 F(\vz(t)) \Big) \,.
\tag{LA3-HRDE}
\end{align*}

\textbf{LA4-GDA}
For LA4-GDA we have
\begin{equation} 
\begin{aligned}
    \vz_{n+1} &= \vz_n + \alpha(\tilde \vz_{n+4} - \vz_n) \\
    &= \vz_n + \alpha\gamma \Big[ 
    - F(\vz_n) 
    - F\big(\vz_n - \gamma F(\vz_n)\big) \\
    &\quad - F\big( \vz_n - \gamma F(\vz_n) 
    - \gamma F(\vz_n - \gamma F(\vz_n)) \big) \\
    &\quad - \underbrace{
        \begin{aligned}
        &F\big( \vz_n - \gamma F(\vz_n) 
        - \gamma F(\vz_n - \gamma F(\vz_n)) \\
        &\quad - F\big( \vz_n - \gamma F(\vz_n) 
        - \gamma F(\vz_n - \gamma F(\vz_n)) \big) \big)
        \end{aligned}
    }_{(\star)}
    \Big]
\end{aligned}
\end{equation}

Similarly, by doing TE in coordinate space for ($\star$),  we get:
\begin{align*}
    F\big( \vz_n - \gamma F(\vz_n) &- \gamma F(\vz_n - \gamma F(\vz_n)) \big) 
    - F\big( \vz_n - \gamma F(\vz_n) -  \gamma F(\vz_n - \gamma F(\vz_n)) \big) \\= 
    &F(\vz_n) + J(\vz_n)(-3\gamma F(\vz_n) + 3\gamma^2 J(\vz_n)F(vz_n)) \\&\quad + \frac{9}{2} \gamma^2 F(\vz_n)^T H(\vz_n) F(\vz_n)
\end{align*}
The final HRDE is as follows 
\begin{align*}
\dot \vz (t) + \frac{1}{2} \gamma \ddot{\vz} (t) + \frac{1}{6} \gamma^2 \dddot{\vz}(t) 
&= \alpha \Big( -4F(\vz(t)) + 6\gamma J(\vz(t))F(\vz(t)) \\
&\quad - \frac{14}{2} \gamma^2 F(\vz(t))^T H(\vz(t))F(\vz(t)) 
      \\& \quad - 4\gamma^2 J(\vz(t))^2 F(\vz(t)) \Big) \,.
\tag{LA4-HRDE}
\end{align*}

\textbf{LA5-GDA}
For LA5-GDA we have
\begin{equation} 
\begin{aligned}
    &\vz_{n+1} = \vz_n + \alpha(\tilde \vz_{n+5} - \vz_n) \\
    &= \vz_n + \alpha\gamma \Big[ 
    - F(\vz_n) 
    - F\big(\vz_n - \gamma F(\vz_n)\big) \\
    &\quad - F\big( \vz_n - \gamma F(\vz_n) 
    - \gamma F(\vz_n - \gamma F(\vz_n)) \big) \\
    &\quad - F\big( \vz_n - \gamma F(\vz_n) 
    - \gamma F(\vz_n - \gamma F(\vz_n)) \\
    &\qquad - F\big( \vz_n - \gamma F(\vz_n) 
    - \gamma F(\vz_n - \gamma F(\vz_n)) \big) \big) \\
    &\quad - \underbrace{
        \begin{aligned}
        &F\big(\vz_n - \gamma F(\vz_n) 
        - \gamma F\big(\vz_n - \gamma F(\vz_n)\big) \\
        &\quad - \gamma F\big( \vz_n - \gamma F(\vz_n) 
        - \gamma F(\vz_n - \gamma F(\vz_n)) \big) \\
        &\quad - \gamma F\big( \vz_n - \gamma F(\vz_n) 
        - \gamma F(\vz_n - \gamma F(\vz_n)) \\
        &\qquad - F\big( \vz_n - \gamma F(\vz_n) 
        - \gamma F(\vz_n - \gamma F(\vz_n)) \big) \big) 
        \big)
        \end{aligned}
    }_{(\star)}
    \Big]
\end{aligned}
\end{equation}

Similarly, by doing TE in coordinate space for ($\star$),  we get:
\begin{align*}
    F\big( \vz_n - \gamma F(\vz_n) - \gamma F(\vz_n - &\gamma F(\vz_n)) \big) 
    - F\big( \vz_n - \gamma F(\vz_n) -  \gamma F(\vz_n - \gamma F(\vz_n)) -\ldots \big) \\ =
    F(\vz_n) + & J(\vz_n)(-4\gamma F(\vz_n) + 5\gamma^2 J(\vz_n)F(vz_n)) 
 \\& + \hfill \frac{16}{2} \gamma^2 F(\vz_n)^T H(\vz_n) F(\vz_n)
\end{align*}

The final HRDE is as follows
\begin{align*}
\dot \vz (t) + \frac{1}{2} \gamma \ddot{\vz} (t) + & \frac{1}{6} \gamma^2 \dddot{\vz}(t) 
= \alpha \Big( -5F(\vz(t)) + 10\gamma J(\vz(t))F(\vz(t)) \\
&\quad - \frac{30}{2} \gamma^2 F(\vz(t))^T H(\vz(t))F(\vz(t)) 
       - 9\gamma^2 J(\vz(t))^2 F(\vz(t)) \Big) \,.
\tag{LA5-HRDE}
\end{align*}

\textbf{LAk -GDA}  The coefficients of the terms for LAk-GDA with $k=2,3,4,5,...$ is given as follows
\renewcommand{\arraystretch}{1.5}
\begin{table}[ht]
\centering
\begin{tabular}{cccc}

$F(\vz_n)$ & $J(\vz_n)F(\vz_n)$ & $F(\vz_n)^T H(\vz_n)F(\vz_n)$  & $(J(\vz_n))^2 F(\vz_n)$\\
\hline
2 & 1 & $\frac{1}{2}$ & 0 \\
3 & 3 & $\frac{5}{2}$ & 1\\
4 & 6 & $\frac{14}{2}$ & 4\\
5 & 10 & $\frac{30}{2}$ & 9 \\
\end{tabular}

\label{coeff-table}
\end{table}

The final $\mathcal{O}(\gamma^2)$ HRDE of LAk-GDA is as follows
\begin{equation}
    \label{eq:LAk-HRDE-second}
    \begin{aligned}
    \dot \vz (t) + \frac{1}{2} \gamma \ddot{\vz} (t) + \frac{1}{6} \gamma^2 \dddot{\vz}(t) 
    &= \alpha \Big( -kF(\vz(t)) + (\sum_{i=0}^{k-1} i) \gamma J(\vz(t))F(\vz(t)) \\
    &\quad - \frac{1}{2}(\sum_{i=0}^{k-1} i^2) \gamma^2 F(\vz(t))^T H(\vz(t))F(\vz(t)) \\
    &\quad - (k-2)^2\gamma^2 J(\vz(t))^2 F(\vz(t)) \Big) \,.
    \end{aligned}
    \tag{LA-$\gamma^2$-HRDE}
\end{equation}
\subsection{\texorpdfstring{Convergence of \ref{eq:LAk-HRDE-second} for simplified \eqref{eq:quadratic_game}}{Convergence of Eq. (X) for Eq. (Y)}}
\label{sec:Convergence_Analysis_ogamma2}
We scale down the \ref{eq:quadratic_game} to the following problem which is a simplified instance of \ref{eq:quadratic_game} ($\vB = \vC = (1-\beta I), \vA = \beta I, \vx =xI, \vy = yI$)
\begin{align*}
    \min_{\substack{x}} \max_{\substack{y}} \ (1-\beta) x^2 + \beta x y - (1-\beta)y^2
\end{align*}
where $\beta \in [0,1]$. 
We follow the same framework for convergence analysis as in Appendix \ref{sec:convergence}. The characteristic equation is as follows
\begin{equation}
\label{eq:polynomial}
    (\frac{1}{6}\gamma^2 s^3 + \frac{1}{2} \gamma s^2 + s + T)^2 + L^2 = 0 /,.
\end{equation}

where $s \in \mathbb{C}$ represents the complex frequency and $T = 2\alpha k(1-\beta) - 2\alpha k(k-1)\gamma (1-\beta)^2 + \alpha \frac{k(k-1)}{2} \gamma \beta^2 +8\alpha (k-2)^2 \gamma^2 (1-\beta)^3 -6\alpha (k-2)^2 \gamma^2 \beta^2 (1-\beta)$ and $L=-\alpha \beta k + 2\alpha k (k-1) \gamma \beta (\beta-1) - \alpha (k-2)^2 \gamma ^2 (4\beta (1-\beta)^2 - \beta^3 + 8\beta (1-\beta)^2$

We compare the polynomial \ref{eq:polynomial} with the following 
\[
P(s) = a_6 s^6 + a_5 s^5 + a_4 s^4 + a_3 s^3 + a_2 s^2 + a_1 s + a_0 /,.
\]
to get the following coefficients 
\begin{align*}
    a_6 = \frac{1}{36} \gamma^4\\
    a_5 = \frac{1}{6}\gamma^3\\
    a_4 = \frac{7}{12} \gamma^2\\
    a_3 = \frac{1}{3}\gamma^2 T + \gamma\\
    a_2 = \gamma T + 1\\
    a_1 = 2T\\
    a_0 = T^2 + L^2
\end{align*}
The corresponding Routh–Hurwitz table is:

\[
\begin{array}{c|cccc}
s^6 & a_6 & a_4 & a_2 & a_0 \\
s^5 & a_5 & a_3 & a_1 & 0 \\
s^4 & b_1 & b_2 & b_3 & 0 \\
s^3 & c_1 & c_2 & c_3 & 0 \\
s^2 & d_1 & d_2 & 0 & 0 \\
s^1 & e_1 & e_2 & 0 & 0 \\
s^0 & f_1 & 0 & 0 & 0 \\
\end{array}
\]

where the elements are computed as:

\[
b_1 = \frac{a_5 a_4 - a_6 a_3}{a_5}, \quad
b_2 = \frac{a_5 a_2 - a_6 a_1}{a_5}, \quad
b_3 = a_0
\]

\[
c_1 = \frac{b_1 a_3 - a_5 b_2}{b_1}, \quad
c_2 = \frac{b_1 a_1 - a_5 b_3}{b_1}
\]

\[
d_1 = \frac{c_1 b_2 - b_1 c_2}{c_1}, \quad
d_2 = a_0
\]

\[
e_1 = \frac{d_1 c_2 - c_1 d_2}{d_1}
\]

\[
f_1 = d_2
\]

For stability, all first-column elements (\( b_1, c_1, d_1, e_1, f_1 \)) must be positive.
i.e.
\begin{equation}
    \label{eq:convergence-o-gamma2}
    \tag{QD-Cond-2}
    b_1>0 \quad \land \quad c_1 >0 \quad \land \quad d_1 >0 \quad \land \quad e_1 >0 \quad \land\quad f_1 >0 \,.
\end{equation}

\section{Settings of Numerical Experiments}
We describe the experimental setting for validating our convergence condition \ref{eq:convergence-condition} and \ref{eq:convergence-inequality-full}, and to empirically show the improvement in the tightness of the convergence bound \ref{eq:convergence-o-gamma2} over \ref{eq:convergence-inequality-full}.  
\subsection{Empirical Verification of \ref{eq:convergence-condition} and \ref{eq:convergence-inequality-full}}
\label{app:experiments-1}

In this section, we describe the setting we used for Figure~\ref{fig:convergence-divergence} in the main part.

We run Lookahead on \ref{eq:bilinear_game} for 200 iterations, using stepsize $\gamma = 0.1$ and $k=5$, with $\alpha = 0.5, 0.8, 0.9$. We observe that the algorithm converges to the Nash Equilibrium steadily for $\alpha =0.5$ which satisfies \ref{eq:convergence-condition}, but diverges for $\alpha = 0.9$ as expected.

To verify \ref{eq:convergence-inequality-full}, we simplify the game \ref{eq:quadratic_game} to a purely potential game, by setting $\vB = \vC =\mathbf{0}_{\frac{d}{2} \times \frac{d}{2}}$. \ref{eq:convergence-inequality-full} reduces to the condition $\gamma < \frac{1}{k-1}$. We run Lookahead on \ref{eq:bilinear_game} for 50 iterations, using $k=5, \alpha =0.5$ and varying $\gamma =0,1, 0.25, 0.5$. We observe that even though $\gamma=0.5, k=5$ does not satisfy \ref{eq:convergence-inequality-full}, Lookahead converges to the equilibrium.

\subsection{Empirical Verification of \ref{eq:convergence-o-gamma2}}
\label{app:experiments-2}

In this section, we describe the setting we used for Figure~\ref{fig:3} in the main part.

In order to empirically confirm that $O(\gamma^2)$ HRDE yields a tighter convergence bound, we scale down the \ref{eq:quadratic_game} to the following problem 
\begin{align*}
    \min_{\substack{\vx}} \max_{\substack{\vy}} \ (1-\beta) \vx^2 + \beta \vx \vy - (1-\beta)\vy^2
\end{align*}
where $\beta \in [0,1]$. 

Varying the parameter $\beta$ from 1 to 0 transitions the game from being purely Hamiltonian to purely potential, effectively modulating the rotational component of the game's dynamics. Let $\gamma$ denote the gradient descent stepsize that satisfies the convergence condition specified in Equation~\ref{eq:convergence-o-gamma2}, and define an alternative stepsize $\gamma' = 1.5 \times \gamma$. Let $d$ represent the distance to the equilibrium. We evaluate the performance of Lookahead combined with gradient descent using both stepsizes $\gamma$ and $\gamma'$, and denote the corresponding distances to equilibrium as $d_{\gamma}$ and $d_{\gamma'}$, respectively. We define the \textit{condition error} as $d_{\gamma} - d_{\gamma'}$, which quantifies the relative effectiveness of the theoretically derived stepsize from the convergence conditions \ref{eq:convergence-condition} and \ref{eq:convergence-inequality-full}. A negative condition error indicates that the stepsize chosen according to the convergence condition leads to faster convergence. We plot this metric across the full range of $\beta$ values for both the $\mathcal{O}(\gamma)$ regime (with the associated convergence condition given in Equation~\ref{eq:convergence-inequality-full}) and the $\mathcal{O}(\gamma^2)$ regime, corresponding to the HRDE condition in Equation~\ref{eq:convergence-o-gamma2}. 

We observe that for $\beta < 0.4$, the condition error is positive for the plot corresponding to the $\mathcal{O}(\gamma)$, while it remains negative for the all values of $\beta$ for the  $\mathcal{O}(\gamma^2)$ curve.

\end{document}